\documentclass[10pt]{article}
\usepackage{amsmath,amssymb,amsthm,graphicx}
\usepackage{url}
\usepackage[numeric]{amsrefs}
\usepackage{hyperref}

\newtheorem{thm}{Theorem}[section]

\newtheorem{lem}{Lemma}[section]
\newtheorem{prop}[thm]{Proposition}

\newtheorem{conj}[thm]{Conjecture}

\newcommand{\thmref}[1]{Theorem~\ref{#1}}
\newcommand{\lemref}[1]{Lemma~\ref{#1}}
\newcommand{\propref}[1]{Proposition~\ref{#1}}

\newcommand{\secref}[1]{Section~\ref{#1}}

\newcommand{\bx}{\hfill$\square$\vspace{.6cm}}

\numberwithin{equation}{section}

\renewcommand\a{\alpha}         
\renewcommand\b{\beta}
\newcommand\g{\gamma}
\renewcommand\d{\delta}
\newcommand\e{\varepsilon}
\renewcommand\l{\lambda}
\renewcommand\L{\Lambda}
\newcommand\D{\Delta}
\newcommand\G{\Gamma}
\newcommand\f{\frac}

\newcommand{\Z}{{\mathbb{Z}}}
\newcommand{\R}{{\mathbb{R}}}
\newcommand{\C}{{\mathbb{C}}}
\newcommand{\A}{{\mathbb{A}}}
\newcommand{\Q}{{\mathbb{Q}}}
\newcommand{\U}{{\mathbb{H}}}

\renewcommand\Re{\mbox{Re~}}
\renewcommand\Im{\mbox{Im~}}

\newcommand{\ttwo}[4]{\left(\begin{array}{cc}
{#1} & {#2} \\ {#3} & {#4} \end{array} \right)}

\renewcommand\i{^{-1}}

\renewcommand\({\left(}         
\renewcommand\){\right)}

\begin{document}

\title{Riemann's Zeta Function and Beyond}
\author{Stephen S. Gelbart\thanks{Partially supported by the
Minerva Foundation.} ~and Stephen D. Miller\thanks{Supported by
NSF grant DMS-0122799.}}
\date{}

\maketitle

\begin{abstract}  In recent years
$L$-functions and their analytic properties have assumed a central
role in number theory and automorphic forms.  In this expository
article, we describe the two major methods for proving the
analytic continuation and functional equations of $L$-functions:
the method of integral representations, and the method of Fourier
expansions of Eisenstein series.  Special attention is paid to
technical properties, such as boundedness in vertical strips;
these are essential in applying the converse theorem, a powerful
tool that uses analytic properties of $L$-functions to establish
cases of Langlands functoriality conjectures.  We conclude by
describing striking recent results which rest upon the analytic
properties of $L$-functions.
\end{abstract}

\begin{flushright}
\emph{Dedicated to Ilya Piatetski-Shapiro, with admiration}
\end{flushright}

 \tableofcontents

\section{Introduction}
\label{Intro}

In 1859 Riemann published his only paper\footnote{See \cite{Riem},
and \cite{Edwards}, for translations.} in number theory, a short
ten-page note which dramatically introduced the use of complex
analysis into the subject. Riemann's main goal was to outline the
eventual proof of the Prime Number Theorem

$$ \pi(x) \ \  = \ \ \#\,\{\hbox{primes }p \le x \} \ \sim \
\f{x}{\log{x}}~,\ \ \ \ ~x\rightarrow\infty \, ,$$
$$\hspace{-1.8 cm}\text{\it{i.e.}}\ \ \ \ \ \ \ \ \ \ \
\lim_{x\rightarrow \infty} \, \pi(x)\ \f{\log{x}}{x} \ \ = \ \ 1
\, ,
$$
 by counting the primes using complex integration (the proof was
  completed half a century later by Hadamard and de la
Vallee Poussin). Along this path he first shows that his
$\zeta$-function, initially defined in the half-plane $\Re(s)>1$
by
 $$\zeta (s) \ \ = \ \  \sum_{n=1}^{\infty}\,\frac{1}{n^{s}}
  \ \ =  \ \ \prod_{p
 \text{~(prime)}}
             \f{1}{1-\f{1}{p^s}} \, ,$$
 has a meromorphic continuation to $\C$.  Secondly, he proposes
what has remained as perhaps the most-famous unsolved problem of
our day:

\vspace{.5 cm}
 {\em {\bf The Riemann Hypothesis:}
 $\zeta(s)\neq 0$ for $\Re{s}>1/2$}.\vspace{.5 cm}
\newline
For more on the history of $\zeta$ and Riemann's work, the reader
may consult \cite{Cartier,Daven,Edwards,Weil-history}.  Our role
here is not so much to focus on the {\em zeroes} of $\zeta(s)$,
but in some sense rather on its {\em poles}.  In particular, our
emphasis will be on explaining how we know that $\zeta(s)$ extends
meromorphically to the entire complex plane, and satisfies the
functional equation
      $$\xi(s)\ :=\ \pi^{-\frac{s}{2}}\, \Gamma
      \, \left(\frac{s}{2}\right)\, \zeta(s)
      \ =\ \xi(1-s)\,.$$
 It is one purpose of this paper to give two separate
treatments of this assertion. We want also to characterize the
$\zeta$-function as satisfying the following three classical
properties (which are simpler to state in terms of $\xi(s)$, the
{\em completed} $\zeta$-function).
\begin{itemize}
\item{{\bf E}ntirety ({\bf E}): $\xi(s)$ has a meromorphic continuation to
  the entire complex plane, with simple poles at $s=0$ and 1.}

  \item {\bf F}unctional {\bf E}quation ({\bf FE}): $\ \xi(s) \ = \ \xi(1-s)$.

\item {\bf B}oundedness in {\bf V}ertical strips ({\bf BV}):
 $\ \xi(s)\,+\,\f 1s \,+\, \f{1}{1-s} \ $ is bounded
in any strip of the form $-\infty<a<\Re(s)<b<\infty$ (i.e.
$\xi(s)$ is bounded in vertical strips away from its two poles).

\end{itemize}
A second purpose is to overview how these treatments and
properties extend to $L$-functions assigned to more general groups
such as $GL(n)$, the group of invertible $n\times n$ matrices (the
function $\zeta(s)$ is attached, we shall see, to $GL(1)$). A
major motivating factor for studying these analytic conditions
(especially the technical {\bf BV}) is that they have become
crucial in applications to the Langlands Functoriality
Conjectures, where they are precisely needed in the ``Converse
Theorem,''  which relates $L$-functions to automorphic forms (see
Theorems~\ref{heckecon}, \ref{weilthm}, and \secref{jpss}). More
to the point, the study and usefulness of $L$-functions has
pervaded many branches of number theory, wherein complex analysis
has become an unexpectedly-powerful tool.  In \secref{la} we
discuss the connections with some of the most dramatic recent
developments, including the modularity of elliptic curves,
progress towards the Ramanujan conjectures, and the results of Kim
and Shahidi.  The two treatments we describe are, in fact, the
major methods used for deriving the analytic properties of
$L$-functions.

\subsection*{The Two Methods}

A first method (\secref{Rie}) of analytic continuation is
Riemann's, initiated in 1859. In fact, it was one of several
different, though similar, proofs known to Riemann; Hamburger, and
later Hecke, moved the theory along remarkably following this line
of attack. Almost a century later, Tate (\secref{Tate}) recast
this method in the modern language of adeles in his celebrated
1950 Ph.D. thesis \cite{Tate}, another famous and important
treatment of $\zeta$-functions.
 The second  -- and
lesser known --  method is via Selberg's ``constant term'' in the
theory of Eisenstein series (\secref{Sel}).  This theory, too, has
an important  expansion: the Langlands-Shahidi method
(\secref{LaSh}).

As we shall see, both  methods take advantage of various (and
sometimes hidden) group structures related to the
$\zeta$-function. They also suggest a wide generalization of the
methods: first to handle Dirichlet $L$-functions,
$\zeta$-functions of number fields, and then quite general
$L$-functions on a wide variety of groups.
 In \secref{modsec} we  begin by explaining this through the
connection between modular forms and $L$-functions.  In fact, this
nexus has been fundamentally important in resolving many classical
problems in number theory.  After surveying the classical theory
of Hecke, we turn to the modern innovations of Langlands.  We have
in mind ideas of Weil, Langlands, Jacquet, Godement,
Piatetski-Shapiro, Shalika, Shahidi, and others.  For broader and
deeper recent reports on the nature of $L$-functions and the
application of their analytic properties, see, for example,
\cite{Iw-Sar,Sar-balt}.

The reader will notice that we have left out many important
properties of the Riemann $\zeta$-function, some related to the
most famous question of all, the Riemann Hypothesis, which can be
naturally restated in terms of $\xi$ as $$\hbox{All $\rho$ such
that $\xi(\rho)=0$ have $\Re\rho=1/2$.}$$ This is because we are
primarily interested in results related to the three properties
{\bf E, BV}, and {\bf FE}. As indicated above, we are also only
following the development of a few approaches (see \cite{Titch}
for many more, though which mainly follow Riemann).  Also, to
lessen the burden on the reader unfamiliar with adeles, we will
more or less describe the historical development in chronological
order, first treating the classical results of Riemann, Hecke,
Selberg, and Weil, before their respective generalizations to
adele groups.

To wit, the paper is organized in three parts.  The first,
Sections 2-4, gives the background on the classical theory:
Section 2 discusses Riemann's theory of the $\zeta$-function and
its analytic properties; Section 3 focuses on Hecke's theory of
modular forms and $L$-functions; and Section 4 centers on
Selberg's theory of non-holomorphic Eisenstein series.  The second
part of the paper redescribes these topics in more modern, adelic
terms. Section 5 leads off with a short introduction to the
adeles. Sections 6, 7, and 8 then give a parallel discussion of
the respective topics of Sections 2, 3, and 4, but in a much more
general context. Finally, the last part of the paper is Section 9,
where we recount some recent results and applications of the
analytic properties of $L$-functions.  Sections 8 and 9 are quite
linked, in that many of  the recent developments and analytic
properties used in section 9 come from the Langlands-Shahidi
method, the topic of Section 8. However, the latter is quite
technical, and we have made an attempt to make Section 9
nonetheless accessible without it.

 A word is in order about what we {\em don't} cover.
Because our theme is the analytic properties of $L$-functions, we
have left out a couple of important and timely topics that lie
somewhat outside our focus. Chiefly among these are some
developments towards the Langlands conjectures, for example the
work of Lafforgue \cite{lafforgue} over function fields.  This is
mainly because the analytic properties of $L$-functions in the
function field setting were long ago established by Grothendieck
(see \cite{katz}), and are of a significantly different nature.
Some resources to learn more about these additional topics include
\cite{MR92j:11045,arthur,MR98d:22017,edfrenkel,
laumon,MR2003c:11051,MR2002e:22024,MR99c:11140,knapp,Rogawski,Bernstein-Gelbart}.

 Before starting, we wish to thank J.
Bernstein, J. Cogdell, A. Cohen, K. Conrad, W. Duke, H. Dym, E.
Lapid, A. Lubotzky, B. Mazur, S.J. Miller, A. Reznikov, B.
Samuels, P. Sarnak, G. Schectman,  F. Shahidi, and the referee for
many very helpful comments.

\section{Riemann's Integral Representation (1859)}
\label{Rie}

As we mentioned in the introduction, Riemann wrote only a single,
ten-page paper in number theory \cite{Riem}.  In it he not only
initiated the study of $\zeta(s)$ as a function of a complex
variable, but also introduced the Riemann Hypothesis and outlined
the eventual proof of the Prime Number Theorem! At the core of
Riemann's paper is the {\em Poisson summation formula}
\begin{equation}
\label{psf}
  \sum_{n \, \in \, \Z} \, f(n)  \ \ = \ \
  \sum_{n\, \in \, \Z}\,\widehat{f}(n) \ ,
\end{equation}
which relates the sum over the integers of a function $f$ and its
Fourier transform
\begin{equation}
\label{fhatdef}
 \widehat{f}(r) \ \ = \ \
  \int_{\R}\,f(x)\,e^{-2\, \pi \, i\,  r\,
  x} \,dx\,.
\end{equation}
The Poisson summation formula is valid for functions $f$ with
suitable regularity properties, such as Schwartz functions: smooth
functions which, along with all their derivatives, decay faster
than any power of $\f{1}{|x|}$ as $|x|\rightarrow\infty.$ However,
by temporarily neglecting such details,
  one can in fact quickly see why the Poisson summation formula
implies the functional equation for $\zeta(s)$, at least on a
formal level. Indeed, let $f(x)\ =\ |x|^{-s}$, so that
\begin{equation}
\label{powerft} \aligned
 \widehat{f}(r) \ \ :=& \ \ \int_\R \, |x|^{-s} \, e^{-2 \, \pi \,i\,
r\, x}\,dx \\
\ \ =& \  \ \ |r|^{s-1}\,G(s)\, ,
\endaligned
\end{equation}
where
\begin{equation}\label{Gdef}
    G(s) \ \ = \ \ \int_\R \,|x|^{-s} \,e^{-2 \,\pi \, i \,
    x}\,dx\,.
\end{equation}
Using the convention that $|0|^s =0$, we can already see from the
Poisson Summation Formula  that
\begin{equation}\label{fakefe}
  2\, \zeta(s)\ \ =\ \ 2 \,G(s)\, \zeta(1-s)\,,
\end{equation}
a functional equation relating $\zeta(s)$ to $\zeta(1-s)$.  In
fact the integral (\ref{Gdef}) is a variant of the classical
$\G$-integral
\begin{equation}\label{gammaintdef}
    \G(s) \ \ = \ \ \int_0^\infty \,e^{-x}\,x^{s-1}\,dx\ ,\ \ \ \
    \ \Re{s}>0
\end{equation}
and can be shown to equal
\begin{equation}\label{Gformula}
    G(s) \ \ = \ \ \f{\pi^{(s-1)/2}\,\G(\f{1-s}{2})}{\pi^{-s/2}\,\G(\f
    s2)}\ ,
\end{equation}
at least in the range $0<\Re{s}< 1$ (see \cite{GR} or
\cite[p.\,73]{Daven}).  Thus the functional equation
(\ref{fakefe}) is formally identical to Riemann's functional
equation
\begin{equation}\label{riemannfe}
    \xi(s) \ \ = \ \
    \pi^{-s/2}\,\G(\textstyle{\f{s}{2}})\,\zeta(s) \ \ = \ \
    \xi(1-s)\, .
\end{equation}

  Of course, neither sum defining $\zeta(s)$ in
 (\ref{fakefe}) converges when the other does, much less in
the range $0<\Re{s}<1$ where we computed $G(s)$.  Indeed,  the
functional equation cannot be proven in the absence of some form
of analytic continuation beyond the region where $\sum_{n=
1}^\infty n^{-s}$ converges.  The argument sketched here for the
functional equation seems to have been first considered by
Eisenstein, who succeeded in proving the functional equation not
for $\zeta(s)$ itself, but for a closely related Dirichlet
$L$-function (for these, see (\ref{dirlfuncdef}) and
\cite{Daven}). Andr\'e Weil has written historical accounts
\cite{Weil-history,Weil-eisen} which suggest that Riemann was
himself  motivated by Eisenstein's papers to analyze $\zeta(s)$ by
Poisson summation. The rigorous details omitted from  the above
formal summation argument
 can be found in \cite[\S 5]{ms-inforder}.

\subsection{Mellin Transforms of Theta Functions}
\label{mellthet}

Riemann's own, rigorous argument proceeds by applying the Poisson
summation formula (\ref{psf}) to the Gaussian $f(x)=e^{-\pi\, x^2
\, t}$, $t>0$, whose Fourier transform is
$$\widehat{f}(r)\ \ = \ \ \f{1}{\sqrt{t}}\ e^{-\,\pi\, r^2/t}\,.$$  The
Gaussian is  a Schwartz function, and can be legitimately inserted
in the Poisson summation formula.  Its specific choice is not
 absolutely essential,
but rather a matter of convenience, as we will see in
\secref{Tate}.  However, it was an inspired selection by Riemann,
in that it is connected to the theory of modular forms (see
\secref{Hecke II}).  Thus Riemann's contribution to the functional
equation went far beyond simply making a formal argument rigorous
-- it launched the link between modular forms and $L$-functions
that remains at the forefront of much mathematical activity a
century and a half later.

By applying Poisson summation to $f(x)=e^{-\pi\,x^2\,t}$    one
thus obtains Jacobi's transformation identity
\begin{equation}\label{jacid}
\theta(i\,t) \ \ =\ \ \f{1}{\sqrt{t}}\
\theta(\textstyle{\frac{i}{t}})\, ,
\end{equation}
where $$\theta(\tau) \ \ = \ \ \f12 \,\sum_{n\, \in \, \Z} e^{\,
\pi\, i\, n^2\,  \tau}  \ \ = \ \ \f12 \ +\ \sum_{n\,=\,1}^\infty
e^{\,\pi\,  i\,  n^2 \, \tau}$$ (more later in \secref{Hecke II}
on $\theta$ as a function of a complex variable for
$\Im{\tau}>0$). Riemann then obtained an integral representation
for $\xi(s)\,=\,\pi^{-s/2}\G(\f s2)\zeta(s)$  as follows:
\begin{equation*}
\aligned\G(s) \ \   & =  \ \  \int_{0}^{\infty} \,t^{s-1} \,e^{-t}
\,
   \,dt \, ,& \Re{s}>0 \\
\pi^{-s}\,\Gamma(s)\,\zeta(2s) \ \
           & = \ \ \sum_{n=1}^{\infty} \
           \int_{0}^{\infty}({\pi} \, {n^2})^{-s}
           \, t^{s-1}\,
 e^{-t}    \,dt\,,& \ \ \ \ \Re{s}>1/2 \\
& = \ \ \int_{0}^{\infty} t^{s-1}\ (\theta(it)-
\textstyle{\frac{1}{2}}) \  dt&
\endaligned
\end{equation*}
\begin{eqnarray*}
&=&\int_{1}^{\infty} t^{\,s-1}\,(\theta(it) -
{\textstyle{\frac{1}{2}}}) \ dt \left.\, \ \ + \  \ \int_{0}^{1}
t^{s-1}\,\theta(it) \, dt  \ \ - \  \
\,\frac{t^s}{2s}\,\right|_{0}^{1}
 \\
&=&\int_{1}^{\infty} t^{s-1}\,(\theta(it) -
{\textstyle{\frac{1}{2}}}) \, dt  \  \ + \  \ \int_{1}^{\infty}
 t^{-s-1}\,\theta(\textstyle{\frac{i}{t}}) \ dt \ \  - \  \  \frac{1}{2s}
  \\
&=&\int_{1}^{\infty} (t^{s-1}\,+\, t^{1/2-s-1})\,(\theta(it)
-{\textstyle{\frac{1}{2}}}) \ dt  \  \ - \  \  \frac{1}{2s} \ \ -
\ \ \frac{1}{1-2s}\,.
\end{eqnarray*}
Indeed, replacing $s$ by $s/2$, the above expression reads
\begin{equation}\label{riemannintegral}
\pi^{-\frac{s}{2}}\,\Gamma({\textstyle \frac{s}{2}})\,\zeta(s)
 \ \ = \ \ \int_{1}^{\infty} (t^{s/2-1}+t^{(1-s)/2-1})
(\theta(it)-\frac{1}{2}) \ dt  \ \ - \ \ \f 1s \ \  - \ \ \f
1{1-s}\, .\end{equation} The integral representation
(\ref{riemannintegral}) allows us to conclude the main analytic
properties mentioned in the introduction:

\begin{thm}\label{zetaebvfe} The function
   $$\xi(s) \ = \ \pi^{-\frac{s}{2}}\, \Gamma({\textstyle \frac{s}{2}}
   )\, \zeta(s)$$
satisfies properties {\bf E, BV}, and {\bf FE} of \secref{Intro}.
\end{thm}

{\bf Proof:}  We first note that
\begin{equation}\label{thetbd}
\theta(it)\,- \,\f 12  \  \ = \  \ \sum_{n=1}^\infty e^{-\pi \,
n^2 \, t} \ \  \le  \ \ \sum_{n=1}^\infty e^{-\pi \, n \, t}  \ \
= \ \ \f{e^{-\pi \, t}}{1-e^{-\pi \, t}} \ \ = \ \ O(e^{-\pi \,
t}) \ \,
\end{equation}
 for $t\ge 1$.\footnote{The notation $A=O(B)$ indicates that there
exists some absolute constant $C>0$ such that $|A|\le C\cdot B$.}
Since
$$\,\int_1^\infty |\, t^{\,s} \, e^{-\pi\, t} \, |\ dt \,  \
 \ \le  \ \
\int_1^\infty t^{\,b} \, e^{-\pi \,  t}\,dt \   \ < \ \  \infty$$
for $\Re{s} \le b$, the integral in (\ref{riemannintegral})
converges -- for {\em any} value of $s$ -- to an entire function
which is bounded for $s$ in vertical strips. Thus $\xi(s)$ is
meromorphic with only simple poles at $s=0$ and 1, demonstrating
properties {\bf E} and {\bf BV}. Having established that
(\ref{riemannintegral}) gives an analytic continuation, we may
conclude that $\xi(s)=\xi(1-s)$ (property {\bf FE}) because of the
symmetry present in (\ref{riemannintegral}). \bx

\subsection{Hecke's Treatment of Number Fields (1916)}
\label{Hecke I}

In this section we shall briefly describe Hecke's generalization
\cite{Hecke} of Riemann's work to certain zeta functions
associated to number fields (that is, finite extensions of $\Q$).
These subsume Riemann's $\zeta$-function, as well as the related
Dirichlet $L$-functions.  The latter are simply Dirichlet series
\begin{equation}\label{dirlfuncdef}
L(s ,  \chi  ) \ \ = \ \ \sum_{n=1}^\infty \, \f{\chi(n)}{n^{s}}\
. \end{equation}
 Here  $\chi$ is a ``Dirichlet character'', meaning a non-trivial function
$\chi:\, \Z \rightarrow \C$ which: (i) is periodic modulo some
integer $N$; (ii) obeys $\chi(nm)=\chi(n)\chi(m)$ (``complete
multiplicativity''); and (iii) vanishes on integers sharing a
common factor with $N$. A Dirichlet character can equally be
thought of as a homomorphism from $(\Z/N\Z)^*$ to $\C^*$, extended
to $\Z$ as a periodic function that  vanishes on $\{n \mid
(n,N)>1\}$.  The Dirichlet $L$-functions $L(s,\chi)$ satisfy the
properties {\bf E, BV,} and {\bf FE} analogous to those of
$\zeta(s)$ (which corresponds to the trivial character); for a
complete discussion and precise analog of \thmref{zetaebvfe}, see
\cite{Daven}.

Our goal here is to describe the generalizations of $\zeta(s)$ and
$L(s,\chi)$ that are the objects of Hecke's work, in some sense
following the earlier exposition in \cite{Gel-langlandssurvey}.
This will necessitate some algebraic background; accordingly this
section requires some familiarity with the concepts involved.
However, it is not essential to the rest of the paper, and readers
may wish to skip directly to \secref{ham}, or instead to consult
\cite{Lant,ramval} for definitions and examples.

To describe Hecke's accomplishment, we need to recall some of the
local and global terminology involved.  Let $F$ be a number field,
and ${\mathcal O}_F$ its ring of integers.  We will refer to a
non-archimedean place $v$ of $F$ as a prime ideal $\mathfrak{P}
\subset {\mathcal O}_F$. A fractional ideal of ${\mathcal O}_F$ is
an ${\mathcal O}_F$-submodule $\mathfrak{U}$ such that
$x\mathfrak{U} \subset {\mathcal O}_F$ for some $x \in F^*$.  All
fractional ideals are invertible (i.e. there exists a
  fractional ideal ${\mathfrak{P^{-1}}}$ such that
 ${\mathfrak{P}}{\mathfrak{P^{-1}}} ={\mathcal O}_F$), and  all fractional
  ideals factor uniquely into products of positive and negative
  powers of prime ideals.  We let $ord_{\mathfrak{P}}(x)$
   denote  the  exponent of $\mathfrak{P}$ occurring  in the unique
factorization of the principal ideal $x {\mathcal O}_F$, and set
  $$|x|_v  \ \ = \ \  |x|_{\mathfrak{P}} \ \ = \ \
  ( N{\mathfrak{P}})^{\,-{ord}_{{\mathfrak{P}}}(x)}\,,$$
where $N{\mathfrak{P}}$ is the number of elements in the finite
field ${{\mathcal O}_F}/\mathfrak{P}$.  Any real embedding
$\sigma:F\rightarrow \R$ of $F$ gives rise to a ``real'' infinite
place via the norm $|x|_v = |\sigma(x)|$; complex places are
defined analogously, and the real and complex places together
comprise the archimedean places of $F$.  For each place of $F$,
the norm $|\cdot|_v$ gives a different completion $F_v$ of $F$.
For example, when $F=\Q$, $F_\infty=\R$ and $F_p=\Q_p$, the
$p$-adic numbers (see \secref{adelesec} for much more on this
theme).

We now come to the generalization of a Dirichlet character to the
number field setting:  a Hecke character (also known as a {\em
Gr\"ossencharacter}). We shall think of one as the product of
family of homomorphisms $\chi_v: F_v^*\rightarrow \C^*$, one for
each place of $F$:
$$\chi(x) \ \ = \ \ \prod_v \, \chi_v(x) \, .$$
Two constraints must be made on the family: firstly that $\chi$ be
trivial on $F^*$, i.e. for any $x \in F \subset F_v^*$
 $$\chi(x) \ \ = \ \ \prod_{ v} \, \chi_v(x) \ \ = \ \ 1 \, ,$$
and secondly that all but a finite number of the $\chi_v$ be {\em
unramified}, i.e. trivial on $\{ x \in F_v^* \mid |x|_v = 1 \}$.
If $v$ is such an unramified place, corresponding to a prime ideal
$\mathfrak{P}$, $\chi(\mathfrak{P})$ is defined as
$\chi_v(\varpi_v)$, where $\varpi_v$ is an element of $F_v$ such
that $|\varpi_v|_v = N{\mathfrak{P}}^{-1}$ (a ``uniformizing
parameter'' for $F_v$).   This definition can of course be
extended to ordinary ideals $\mathfrak{U}$ of ${\mathcal O}_F$,
provided they are products of prime ideals corresponding to places
where $\chi_v$ is unramified.   Hecke's (abelian) $L$-series for
the character $\chi$ is then defined as the Dirichlet series
\begin{equation}
\label{heckeslseries}
 L(s,\chi) \ \ = \ \ \sum {\frac{\chi(\mathfrak{U})}{(N\mathfrak{U})^{s}}}
\    \ = \  \ \prod_{\mathfrak{P}} \left( \, 1 \ - \
{\chi}(\mathfrak{P})\,
       (N{\mathfrak{P}})^{-s} \, \right)^{-1}.
       \end{equation}
Here ${\mathfrak{U}}$ is summed over these ordinary ideals of
${\mathcal O}_F$ just mentioned, and the product is only over the
prime ideals corresponding to these unramified places.

When $\chi$ is the trivial character, i.e., $\chi_v=1$ for all
$v$, then $L(s,\chi)$ specializes to be the \emph {Dedekind
zeta-function} $\sum (N\mathfrak{U})^{-s}$  of $F$. For $F=\Q$
this reduces to $\zeta(s)$,  and if  $\chi$ is instead of finite
order, $L_F(s,\chi)$  becomes the Dirichlet $L$-function
(\ref{dirlfuncdef}). Using  very clever and intricate arguments,
Hecke was able to express his $L$-series in terms of generalized
``$\theta$-functions'', and to then derive their analytic
continuation, functional equation, and boundedness in vertical
strips, {\em a la} Riemann.

\subsection{Hamburger's Converse Theorem (1921)}
\label{ham} Now let us return to the Riemann $\zeta$-function. The
next point of the theory is that the {\bf F}unctional {\bf
E}quation for $\zeta(s)$ nearly characterizes it. Indeed,
Hamburger \cite{Ham} showed in 1921 that any Dirichlet series
satisfying $\zeta$'s functional equation \emph{and suitable
regularity conditions} is necessarily a constant multiple of
$\zeta(s)$.  We state Hamburger's \thmref{hamthm} at the end of
this section, but first begin by describing these conditions,
which are closely related to {\bf B}oundedness in {\bf V}ertical
strips.  In fact, our main motivation in describing Hamburger's
theorem here is to explain the role of {\bf BV} and the related
``finite order'' conditions; in the modern picture, these are
crucial for applications  involving the Converse Theorem (see
Sections~\ref{Hecke II}, \ref{Weil}, and \ref{jpss}).

Recall that property {\bf BV} was stated earlier in terms of the
function $\xi(s)=\pi^{-s/2}\G(\f s2)\zeta(s)$.  One may ask how
the individual factors themselves behave as
$|\Im{s}|\rightarrow\infty$. Clearly
\begin{equation}\label{piasym}
  |\pi^{-s/2}| \ \ = \ \ \pi^{-\text{Re}(s)/2} \, ,
\end{equation}
while Stirling's formula states that
\begin{equation}\label{stirup}
\gathered
  |\G(\sigma+it)| \  \ \sim \  \ \sqrt{2\pi}\, |t|^{\sigma-1/2} \,
  e^{-\pi\,
  |t|/2}\,,  \\ \hbox{uniformly for }a \, \le  \, \sigma  \, \le \,
   b, \ \ \ |t| \rightarrow \infty.
\endgathered
\end{equation}
Yet the size of $|\zeta(s)|$ in the critical strip is quite
difficult to pin down.  In fact, one of the central unsolved
problems in analytic number theory is the following and its
generalizations. \vspace{.6 cm}
\\
{\em {\bf  The Lindel\"of Hypothesis:} For any fixed $\e>0$ and
$\sigma\ge 1/2$,
\begin{equation}
\label{lindl} \,\zeta(\sigma\, +\, i\,t)\, \ \ = \ \ O( |t|^{\e})\
\ \ \ \ \hbox{ as }  \ \ \ \  \ |t|\rightarrow \infty.
\end{equation}}  The implied constant in the $O$-notation here depends
implicitly on the value of $\e$. In particular,
$|\zeta(1/2+it)|=O(|t|^\e)$ for $|t|$ large (this case turns out
to be equivalent to (\ref{lindl}) via the Phragmen-Lindel\"of
Principle, \propref{phrag}). The Lindel\"of Hypothesis is implied
by the Riemann Hypothesis, and conversely implies that very few
zeros disobey it (see \cite[\S 13]{Titch}).

Note that by (\ref{piasym}), (\ref{stirup}), and the {\bf
F}unctional {\bf E}quation, the behavior for $\Re(s)\le 1/2$ is
given by
\begin{equation}\label{fnbd}
|\zeta(\sigma+it)| \ \sim  \ |\zeta(1-\sigma-it)|
\left|\f{t}{2\pi}\right|^{1/2-\sigma}, \ \ \ \sigma \text{ fixed},
~~~|t| \text{ large.} \end{equation} The Lindel\"of conjecture is
far out of reach, but we can easily prove (weaker) polynomial
bounds.

\begin{prop}\label{trivzetabd}
$\zeta(s)-\f{1}{s-1}=O(|s|)$ for $\Re{s}\ge 1/2$.
\end{prop}
{\bf Proof:}

For $\Re{s}>1$,
\begin{eqnarray}
\zeta(s) \ - \ \f{1}{s-1} &\ = \ \ & \sum_{n=1}^\infty n^{-s} \ -
\ \int_1^\infty
x^{-s}dx \\
&\ = \ \ &\sum_{n=1}^\infty \int_n^{n+1} (n^{-s}-x^{-s})\,dx\, .
\label{firstline}
\end{eqnarray}
The integrand in (\ref{firstline}) is bounded by
$$|n^{-s}-x^{-s}| \ \ =  \ \
 \left|\, \int_n^x s\,t^{-s-1}\, dt \, \right|
\ \  \le  \ \  |s|\ n^{-\,\scriptstyle{\text{Re}}{\,s}\,-\,1}\,.$$
We conclude that
$$\left|\, \zeta(s)\, - \, \f{1}{s-1}\, \right|  \ \ \le
\ \  \  |s|\ \zeta(\Re s+1)\, ,$$ and so (\ref{firstline}) gives
an analytic continuation of $\zeta(s)-\f{1}{s-1}$ to the region
$\Re s > 0$.  In particular, $|\zeta(s)-\f{1}{s-1}| \,\le \,|s|\,
\zeta(3/2)$ for $\Re{s}\ge 1/2$. \bx

 {\bf Definition:} An entire function $f(s)$ is {\em of order
$\rho$} if
\begin{equation}\label{entdef}
f(s) \ \  = \ \ O({e}^{{|s|}^{\rho+{\epsilon}}}) \ \ \ \ \hbox{
for any $\epsilon>0$.}
\end{equation}
It will turn out that the $\zeta$-function and (conjecturally) all
$L$-functions connected to automorphic forms have order 1.
However, many other generalizations of zeta functions (such as
Selberg's Zeta functions) in fact have order greater than 1.

\begin{prop} The function
  $$ s \, (s-1) \, \pi^{-\frac{s}{2}}\,
  \Gamma({\textstyle \frac{s}{2}} ) \, \zeta(s)$$
is (entire and) of order $1$.
\end{prop}

{\bf Proof:}  By the {\bf F}unctional {\bf E}quation, it suffices
to consider $\Re{s}\ge 1/2$. We have already seen this function is
{\bf E}ntire in \thmref{zetaebvfe}. Another form of Stirling's
Formula gives that

\begin{equation}\label{bigstir}
  \G(s) \ \sim  \ \sqrt{2\pi} \, e^{-s}\,s^{\,s-\f 12}\ = \
  \sqrt{2\pi}\,e^{-s+(s-\f 12)\log{s}}\, , \ \ \  \ \Re{s}\ge 1/2,\
  |s|\rightarrow \infty,
\end{equation} and hence $\G(s)=O(e^{\,M |s|\log|s|})$ for some
$M>0$.  Thus by (\ref{piasym}), (\ref{bigstir}), and
\propref{trivzetabd}
$$s \, (s-1)\, \pi^{-s/2}\, \G({\textstyle \f s2})\,\zeta(s)\
=\ O(|s|^3 e^{\f M2 |s|\log|s|}).$$  Since for any $\e>0$,
$\f{\log|s|}{|s|^{\e}}\rightarrow 0$ as $|s|\rightarrow \infty$,
we conclude
$$s\,(s-1)\,\pi^{-s/2}\,\G({\textstyle \f s2})\,\zeta(s)
\ = \ O(e^{|s|^{1+\e}}).$$\bx

We note that $s\,(s-1)\,\xi(s)$ is not of any  order $\rho<1$, as
can be seen from (\ref{piasym}) and (\ref{bigstir}) as
$s\rightarrow\infty$ along the real numbers -- where $\zeta(s)$ is
always greater than 1.

\begin{thm}\label{hamthm}(Hamburger's Converse Theorem\footnote{Actually
Hamburger proved a more general statement, allowing for an
arbitrary, finite number of poles (see
\cite{Ham},\cite[p.31]{Titch}).})

 Let
 $h(s)\, = \, \sum_{n=1}^{\infty}\ a_n \,n^{-s}$ and
 $g(s)\, = \, \sum_{n=1}^\infty b_n \, n^{-s}$
  be absolutely
convergent for $\Re s >1$, and suppose that both $(s-1)h(s)$ and
$(s-1)g(s)$ are entire functions of finite order. Assume the
functional equation
\begin{equation}\label{hamfe}
\pi^{-\frac{s}{2}}\, \Gamma({\textstyle \frac{s}{2}}) \, h(s) \ \
= \ \  \pi^{-\frac{1-s}{2}} \, \Gamma({\textstyle
\frac{1-s}{2}})\, g(1-s)\,.
\end{equation}
 Then in fact
$h(s)  \,  =  \,  g(s)  \, = \,   a_1\,\zeta(s)$.
\end{thm}

This is the theorem which says that  $\zeta(s)$ is \emph{uniquely}
determined by its functional equation (subject to certain
regularity conditions). Hamburger's theorem was greatly
generalized and enlightened  by Hecke approximately 15 years
later. We will in fact later show how to derive \thmref{hamthm}
from Hecke's method (see the discussion after \thmref{heckecon}.)
See also \cite{psrag,rag}.

The original proof of Hamburger's Theorem relies on the
\emph{Mellin transform} and  {\em inversion formulas}; that is, if
\begin{equation}\label{MT}
    \pi^{-s}\,\Gamma(s)\,\zeta(2s) \ \ = \ \
    \int_{0}^{\infty} t^{s-1}\,
  (\theta(it)-\frac{1}{2}) \, dt \, ,
\end{equation}
then
\begin{equation}\label{MTinv}
    \theta(it)\,-\,\frac{1}{2} \ \ = \ \
    \frac{1}{2{\pi}i} \, \int_{Re(s)\, = \, c}
         t^{-s}\,\( \, \pi^{-s}\, \Gamma(s)\, \zeta(2s)\, \) \, ds
\end{equation}
for sufficiently large $c>0$.  Using the Phragmen-Lindel\"of
principle (\propref{phrag}), \emph {plus} the regularity
conditions of $g(s)$ and $h(s)$, one can fairly directly show that
every $a_k$ is equal to $a_1$; that is, $h(s)=a_1{\zeta(s)}$. By
the way, it is of course known that $(s-1)\zeta(s)$ is entire and
of order 1. This is because both $s(s-1)\pi^{-s/2}\G(s/2)\zeta(s)$
and $\f{1}{s\G(s/2)}$ are entire and of order 1
(\thmref{zetaebvfe}).

\subsection{The Phragmen-Lindel\"of Principle and Convexity Bounds}
\label{frag}

A standard fact from complex analysis, the {\em
Phragmen-Lindel\"of Principle}, can be used to obtain estimates on
$\zeta(s)$ in vertical strips from ones on their edges:

\begin{prop}\label{phrag}(Phragmen-Lindel\"of).  Let $f(s)$ be
meromorphic in the strip $U=\{s \, | \, a\le \Re(s) \le b\}$,
$a,b\in\R$, with at most finitely-many poles.  Suppose that $f(s)$
satisfies the finite order inequality $$f(s) \ \ = \ \
O\(e^{|s|^A}\)\, , \ \ \ \hbox{for some~}A
> 0 \, ,
$$ on $U$ for $|\Im s|$ large, and obeys the estimate
$$f(\sigma+it) \ \ = \ \ O(|t|^M) \ \ \ \text{ for }~\Re s
=a,b, \ \
|t|\rightarrow \infty.$$
 Then
$$f(\sigma+it) \ \ = \ \ O(|t|^M) \ \ \ \text{ for }~a\le \Re{s}\le
 b, \ \ |t|\rightarrow
\infty$$ as well.
\end{prop}

See \cite{Laca} for a detailed exposition and proof of
\propref{phrag}.  An immediate application of the
Phragmen-Lindel\"of Principle is to provide  a standard bound for
$\zeta(s)$ and other $L$-functions in the critical strip.  As an
example, let us note the following bound towards the Lindel\"of
conjecture:

\begin{lem}\label{convex}
For any $\e>0$,
\begin{equation}\label{coneq}\zeta(1/2\,+\,i\,t) \ \
 = \ \ O_\e(t^{1/4+\e})\ ,
~~~ \ \ \ \ \ |t|\rightarrow\infty
\end{equation} where the implied constant depends on $\e$.
\end{lem}
Note that this a sizeable improvement on the trivial bound in
\propref{trivzetabd} towards (\ref{lindl}).

{\bf Proof of \lemref{convex}:}  First we observe that
$$|\, \zeta(1+\e+it) \, | \ \ \le \ \  \sum_{n=1}^\infty |
\,n^{-1\,-\,\e\,-i\,t}\,| \ \ = \ \  \zeta(1+\e)\, ,$$ which is a
positive constant.  By (\ref{fnbd}), which comes from the
functional equation,
$$| \, \zeta(-\e-it) \, |\ \ = \ \ O_\e(|t|^{1/2+\e})\, , \ \ \ \ \ \
|t|\rightarrow\infty \, .$$  Now, set $f(s)=\zeta(s)\zeta(1-s)$,
$a=-\e$, $b=1+\e$, and $M=1/2+\e$.  Because of  the discussion at
the very end of \secref{ham}, the conditions of \propref{phrag}
are met;
 we conclude $|\zeta(1/2+it)\zeta(1/2-it)|=O_\e(|t|^{1/2+\e})$ as
 $|t|\rightarrow\infty$.
   To finish the proof we replace $\e$ by $2\e$, and
observe that $\zeta(\bar{s})=\overline{\zeta(s)}$ because of the
Schwartz reflection principle ($\zeta(s)=\sum n^{-s}$ is real for
$s>1$). \bx

The  estimate (\ref{coneq})  for $\zeta(s)$ has been improved many
times over; however, for general $L$-functions, the bounds given
by the above argument are usually the best  known.  Because
(\ref{coneq}) interpolates between the bounds at $\Re{s}=-1-\e$
and $\e$, results given by this argument are known as the {\em
convexity bounds} for $L$-functions.  A very important problem is
to improve these by {\em breaking convexity}; even slight
improvements to the convexity bounds for more general
$L$-functions -- still falling far short of Lindel\"of's
conjecture -- have had many profound applications. Let's consider,
for example, the possible ways of writing a positive integer as
the sum of three squares. Gauss' famous condition asserts that the
equation
\begin{equation}\label{3squares}
x^2\ + \ y^2\ + \ z^2\ \ = \ \ n\end{equation} is solvable by some
$(x,y,z)\in\Z^3$ if and only if $n$ is not of the form $4^a(8b+7)$
for some integers $a,b\ge 0$
 (see, for example, \cite{Serre}).  Linnik conjectured that the
solutions to (\ref{3squares}) are randomly distributed in the
sense that the sets
\begin{equation}\label{Dn}
  {\mathcal D}_n  \ \ = \ \  \left\{\left. \f{(x,y,z)}{\sqrt{n}}\,
   \right|\  x^2+y^2+z^2\,=\,n \ ,
  \ \ \  ~~~x,\,y,\, z \, \in \, \Z \right\}
\end{equation}
become equidistributed in the sphere $S^2\subset \R^3$ as $n\neq
4^a(8b+7)$ increases. This was in fact proven by W. Duke (see
\cite{Duke,I1,Duke-Rankin}), and can be shown to follow quite
directly from subconvexity estimates on automorphic $L$-functions
(\cite{DFI}), although this was not Duke's original argument. For
a survey of recent results on sub-convexity bounds, see
\cite{Iw-Sar}.

The proof of \lemref{convex} shows the strength of the
finite-order condition.  For it allows us to conclude that
$\xi(s)$ decays rapidly as $|\Im s|\rightarrow\infty$ (and
uniformly so in vertical strips), given only the functional
equation and the absolute convergence of $\zeta(s)$ for $\Re\!(s)$
large. This will be useful in the proofs of Theorems~\ref{hamthm}
and \ref{heckecon}. To wrap up this section, let's formally state
this for future use.

\begin{lem}\label{fogetsbv}
Assume the conditions  of \thmref{hamthm} (notably {\bf E}ntirety,
{\bf F}unctional {\bf E}quation, and the finite order hypothesis).
Then both sides of (\ref{hamfe}) are {\bf B}ounded in {\bf
V}ertical strips.
\end{lem}

We remark that the conclusion of the Lemma does not depend
particularly on the exact form of the {\bf F}unctional {\bf
E}quation (\ref{hamfe}); similar conclusions follow when the
functional equation involves different configurations of
$\G$-functions and powers of $\pi$.

{\bf Proof:}  The assumption of absolute convergence implies that
$$|h(\sigma+it)| \ \ \le \ \  \sum_{n=1}^\infty
|a_n|\,n^{-\sigma} \ \,  < \ \ \infty  \  ,~~~\ \ \ \ \ \
\sigma>1.$$ Then for any $\e>0$, $|h(s)|$ is uniformly bounded in
the range $\Re{s}\ge 1+\e$, as is $|g(s)|$ by symmetry. Using the
{\bf F}unctional {\bf E}quation, we see that both
$$|h(\sigma+it)| \ , \ \ |g(\sigma+it)| \ \ = \ \ O(|t|^{1/2-\sigma})\ ,~~~
\ \ \ \ \ \ |t|\rightarrow\infty$$
 for
$\sigma<-\e$, and uniformly so in vertical strips (see
(\ref{fnbd})).

We  are assuming that $(s-1)g(s)$ and $(s-1)h(s)$ are of finite
order, so  the Phragmen-Lindel\"of Principle (\propref{phrag})
applies. This shows that
$$|g(\sigma+it)| \ , \ \ |h(\sigma+it)|
\ \ = \ \ O(|t|^{1/2+\e}) \, , \ \ \ \ \ \ ~~~\hbox{for }
-\e<\sigma<1+\e.$$ Thus we have shown that in any vertical strip
$a\le \Re{s} \le b$, both $g(s)$ and $h(s)$ are bounded by
$|\Im{s}|^M$ for some $M>0$, as $|\Im{s}|\rightarrow\infty$.
Stirling's estimate (\ref{stirup})  shows that (\ref{hamfe})
decays rapidly as $|\Im{s}|\rightarrow\infty$ in the strip $a\le
\Re{s} \le b$, and  hence is  bounded there. \bx

\section{Modular Forms and the Converse Theorem}\label{modsec}

\subsection{Hecke (1936)}
\label{Hecke II}

As already suggested, Hamburger's Converse Theorem did not become
completely understood until greatly generalized by Hecke in 1936
(\cite[paper \#33]{Hecke}, \cite{Heckenotes}); to describe it, we
thus encounter the notion of the space of modular forms to which
functions like $\theta$ belong. Note that

$$\theta(\tau) \ \ = \ \ \frac{1}{2}\, \sum_{n\,= \, -\infty}^{\infty}
e^{\, {\pi}\, i\, {{n}^2}\, \tau}$$ is holomorphic in the upper
half plane $\Im\!(\tau)>0$; moreover, because it satisfies
(\ref{jacid}) (when $\Re{\tau}=0$), clearly
\begin{equation}\label{thetmod}
\theta\(\f{-1}{\tau}\) \ \ = \ \ \(\f{\tau}{i}\)^{1/2} \,
\theta(\tau) \ , \ \ \ \ \ \theta(\tau + 2) \ \ = \ \
\theta(\tau)\, .\end{equation} These two equations say that
$\theta(\tau)$ is a \emph{modular form of weight $\frac{1}{2}$}
for the group generated by $\tau \mapsto \tau + 2$ and $\tau
\mapsto -\frac{1}{\tau}$.  More generally, a \emph{modular form of
weight $k>0$} and {\em multiplier condition $C$} for the group of
substitutions generated by $\tau \mapsto \tau
 +\lambda$ and $\tau \mapsto -{\frac{1}{\tau}}$ is a holomorphic
function $f(\tau)$ on the upper half plane satisfying
\begin{itemize}
\item      (i) $f(\tau + \lambda) \ = \ f(\tau)$,

\item  (ii) $f(-\frac{1}{\tau}) \ = \ C\,({\frac{\tau}{i}})^k \,
f(\tau)$, and

  \item    (iii) $f(\tau)$ has a Taylor expansion in
$e^{\frac{2\, {\pi} \, i \, {\tau}}{\lambda}}$ (cf. (i)): $f(\tau)
= \sum_{n=0}^{\infty} a_n\, {e^{\frac{\, 2\, {\pi} \, i \,  n  \,
{\tau}}{\lambda}}}$, i.e., $f$ is ``holomorphic at $\infty$''.
\end{itemize}
 We denote the space of such $f$ by
$M(\lambda, k, C)$; $f$ is a {\em cusp form} if $a_0=0$.  For
example, the space $M(2,\f 12, 1)$ is one dimensional, and
consists of multiples of the $\theta$-function.

Now, given a sequence of complex numbers $a_0, a_1, a_2,\ldots$
with $a_n= O(n^d)$ for some $d>0$, and given $\lambda > 0, k>0,
C=\pm 1$, set
\begin{equation}
     \phi(s) \ \ = \ \ \sum _{n=1}^{\infty} \ \frac{a_n}{n^{s}}\, ,
     \end{equation}
     \begin{equation}
     \Phi(s) \ \ = \ \ \(\frac{2{\pi}}{\lambda}\)^{-s} \,  \Gamma(s)
     \,\phi(s) \, ,
     \end{equation}
and
     \begin{equation}\label{fdef}f(\tau) \ \ = \ \
     \sum _{n=0}^{\infty} \,
     a_n\,
{e^{\frac{\,2\,{\pi}\, i  \, n \,
{\tau}}{\lambda}}}.\end{equation} (The $O$-condition on the $a_n$
ensures that $\phi(s)$ converges for $\Re s >d+1$, and that $f(s)$
is holomorphic in the upper half plane. In fact,
$f(\tau)-a_0=O(e^{-\f{2\pi}{\l}\scriptstyle{\text{ Im~}}{\tau}})$
-- see (\ref{thetbd}).)

\begin{thm}(Hecke's Converse Theorem)\label{heckecon}
The following two conditions are equivalent:

(A) ${\Phi(s)} \, + \, \frac{a_0}{s} \, + \, \frac{{C}a_0}{k-s}$
is an entire function which is bounded in vertical strips { \bf
(EBV)}, and satisfies $\Phi(s) \ = \ C \, {\Phi(k-s)}$ {\bf(FE)};

(B) $f$ belongs to $M(\lambda, k, C)$.
\end{thm}

We will come to the proof \thmref{heckecon} shortly, but first
wish to explain the connection to the results of Riemann and
Hamburger.  Riemann's \thmref{zetaebvfe} is an example of the
direction (B)$\Rightarrow$(A). In the other direction, set
$$\gathered \phi(s) \ \ = \ \ \zeta(2s)
 \  \ = \ \ \sum_{n\ge
1}\ (n^2)^{-s} \, , \\  \  C=1 \, , \  k=1/2 \,   ,  \ a_0=1/2 \,
, \, \ \text{and }\, \lambda=2 \,. \endgathered $$ In this special
case, the direction (A)$\Rightarrow$(B) of  \thmref{heckecon}
asserts that $\theta(\tau)$ obeys the modular relations
(\ref{thetmod}). Similarly, \thmref{hamthm} can be derived from
this direction of \thmref{heckecon} as well.  For simplicity,
suppose that the coefficients $a_n$ and $b_n$ in the statement are
equal (these are not the same $a_n$ involved here).  Then
assumptions of \thmref{hamthm} actually match the properties of
$\zeta$ and $\xi$ needed in (A). They guarantee that $\Phi(\f
s2)=\pi^{-s/2}\G(\f s2)h(s)$ is holomorphic in $\Re{s}>0$, except
perhaps for a simple pole at $s=1$.  By the functional equation
(\ref{hamfe}), $\Phi(s)$ has an analytic continuation to $\C$
except for potential simple poles at $s=0$ and $1/2$. Because of
(\ref{hamfe}) the residues of $\Phi(s)$ at those points are
negatives of each other, and thus
${\Phi(s)}+\frac{a_0}{s}+\frac{a_0}{k-s}$ is {\bf E}ntire, where
$a_0$ is the residue of $\Psi(s)$ at $s=k$. \lemref{fogetsbv}
shows the finite order assumption implies that $\Phi(s)$ satisfies
the {\bf BV} condition of (A) as well.
 \thmref{heckecon} therefore produces a
 modular form $f$ in $M(2,1/2,1)$, which is  a
one-dimensional space spanned by $\theta(\tau)$. So $f$ must
 in fact be a multiple of the $\theta$-function, and we  conclude
 that the original Dirichlet series in \thmref{hamthm} are
 multiples of $\zeta$.

{\bf Proof of \thmref{heckecon}:} As in Hamburger's proof of
\thmref{hamthm}, the proof begins by using Mellin inversion (see
(\ref{MTinv})):
\begin{equation}\label{hamstart}
    f(ix)  \, - \,  a_0 \ \  = \ \ \frac{1}{{2}{\pi}i} \,
    \int_{{\sigma} \, = \, c} \ x^{-s} \, \Phi(s) \, ds \,  ,
\end{equation}
for $x>0$, where ${\sigma}=\Re(s)$, and $c$ is chosen large enough
to be in the domain of absolute convergence of $\phi(s)$ (since we
are assuming that $a_n=O(n^d)$, we may take any $c>d+1$).

Assume now (A).
 We first want to first argue that we can
push the line of integration to the left, past $\sigma = 0$,
picking up residues of $ C\,a_0\,x^{-k}$ at $s=k\le c$ and $-a_0$
at $s=0$:
\begin{equation}\label{ham2}
  f(ix) \ - \ C \,{a_0}\,{x^{-k}} \ \   =
   \ \ \frac{1}{{2}\,{\pi}\,i} \,
   \int_{{\sigma} \, = \, k \, - \, c \, < \,0}
                    \ x^{-s} \, \Phi(s) \, ds \, .
\end{equation}

To see this, we need to show that the integral of $\Phi(s)$ over
the horizontal paths $[k-c\pm ir, c\pm ir]$  tend to zero as $r
\rightarrow \infty$. We shall use the {\bf B}oundedness in {\bf
V}ertical strips assumption to prove the integrand decays rapidly
there; in fact, this contour shift here is the primary importance
of the {\bf BV} property. The assumption that $a_n=O(n^d)$ implies
that $\phi(s)$ converges absolutely for $\Re{s} \ge c > d+1$,
where
\begin{equation}\label{phiest}
  |\phi(s)| \ \ \le \ \  \sum_{n=1}^\infty \,
|a_n| \, n^{-c} \ \ = \ \ O(1) \, .
\end{equation}
Stirling's asymptotics (\ref{bigstir}) show that $\Phi(s)$
satisfies the order-one estimate $O(e^{|s|^{1+\e}})$
in~the~region~$\Re{s}\ge c$. By the functional equation, $\Phi(s)$
does as well in the reflected region $\Re{s}\le k-c$, and the {\bf
BV} assumption from (A) handles the missing strip: therefore
$\Phi(s)+\f{a_0}{s}+\f{Ca_0}{k-s}$ is of order one  on $\C$.
 Since
$\f{1}{s\G(s)}$ is entire and of order 1,
$(s-k)\phi(s)=(s-k)(\f{2\pi}{\l})^s\G(s)\i \Phi(s)$ is also entire
and of order 1  (cf. the end of \secref{ham}). The functional
equation
$$\phi(s) \ \ = \ \
C\(\f{2 \, \pi}\l\)^{2s-k}\f{\G(k-s)}{\G(s)} \, \phi(k-s)$$ shows
that
\begin{equation}\label{phiest2}
  |\phi(\sigma+it)| \ \ = \ \ \ O(t^{2c-k}) \, , \ \ \ \sigma
  \, = \, k \, - \, c \, < \, 0
\end{equation}
just as in the proof of \lemref{convex}.  We conclude from the
Phragmen-Lindel\"of Principle (\propref{phrag}) that $\phi(s)$ is
$O(|\Im{s}|^K)$ for some $K$, uniformly as $|\Im{s}|\rightarrow
\infty$ in the strip $k-c \le \Re{s} \le c$.  Since this growth is
at most polynomial, the exponential decay from Stirling's formula
(\ref{stirup})  gives us that $\Phi(\sigma+it)$ decays faster than
any polynomial in $|t|$ as $|t|\rightarrow \infty$, uniformly for
$\sigma$ in the interval $[k-c,c]$. Thus the integrals
$$\int_{k-c+ir}^{c+ir}  \, x^{-s} \, \Phi(s) \, ds \ , \
\int_{k-c-ir}^{c-ir} \, x^{-s} \, \Phi(s) \, ds \ \ \
\longrightarrow \,  0 \hbox{ , as }r\rightarrow \infty \, ,$$ and
the contour shift between (\ref{hamstart}) and (\ref{ham2}) is
valid.

Now, let us resume from (\ref{ham2}) and apply the functional
equation from (A):
\begin{equation*}
\aligned
 f(ix) \ - \ C \,a_o \, x^{-k} \ \ & = &
                \frac{C}{{2}{\pi}i} &\int_{{\sigma} \, = \, k \, - \,
                c \, < \, 0}
                   \ x^{-s} \, \Phi(k-s) \, ds  & \\
   \ \ & = &
       \frac{C}{{2}{\pi}i}& \int_{{\sigma} \,  =  \, c \, > \, k}
       \ x^{s-k} \, \Phi(s) \, ds &
      (\hbox{upon }s\mapsto k-s)  \qquad\qquad \\
 \ \ & = &  C\,{x^{-k}} &\, (f(\frac{i}{x})\, -\,
          {a_0})
         &  \hbox{ by (\ref{hamstart})}\, , \qquad\qquad
\endaligned
\end{equation*}
           or
           $$ f(ix) \ \ = \ \ C\, {x^{-k}}\, f(\frac{i}{x}) \, ,$$
which is property (ii) of the definition of $M(\l,k,C)$.
Properties (i) and (iii) are immediate from the definition of
$f(\tau)$ in (\ref{fdef}), and we conclude the proof that (B)
follows from (A).

Now suppose (B).  We will essentially follow Riemann's original
argument from \secref{mellthet}, using the integral representation
(cf. (\ref{MT}))
$$\Phi(s) \ \ = \ \ \int_{0}^{\infty} \ t^{s-1}
\, (f(it) \, -  \, a_0) \, dt
\, .$$ Then

$$\int_{0}^{1} \ t^{s-1} \, (f(it) \, - \, a_0) \, dt \ \
= \ \    \int_{1}^{\infty}
 t^{-s-1} \,  f(\frac{i}{t}) \, dt
  \ \,  -  \, \
 \left. {a_0} \, \frac{t^s}{s}\, \right|_{0}^{1}
  \ \ \ \ $$

  $$ \ \ \ \ \ \ \ \  \ = \ \  C { \int_{1}^{\infty}}
 t^{k-s-1} \, (f(it)\,-\, a_0) \, dt  \ - \ \f{a_0}{s}  \ - \  \frac{C\, {a_0}}{k-s}\,. $$
Thus

 $$\Phi(s) \   + \  \frac{a_0}{s}  \ + \ \frac{C\,{a_0}}{k-s} \ \ =
 \ \  \int_{1}^{\infty} [ t^{s-1}\,(f(it) \, - \,  a_0)
  \  + \ t^{k-s-1} \, C(f(it)\, - \, a_0)] \ dt\,.$$
This expression is clearly {\bf EBV}, and
$\Phi(s)\,=\,C\,\Phi(k-s)$ (just as in the proof of
\thmref{zetaebvfe}), whence (A).\bx

By reducing a question about Dirichlet series to one about
 modular forms, \thmref{heckecon} represents a great
step forward from Riemann's treatment of $\zeta$.  In particular,
it puts his original argument into a very useful and fruitful
context. Note that a specified type of Dirichlet series is
connected to any modular form satisfying
$$f\(\f{a\tau+b}{c\tau+d}\)=(c\tau+d)^kf(\tau)$$ for

\begin{equation}\label{sl2zdef}
 \ttwo abcd \ \ \in \ \    SL(2,\Z) \ \ = \ \
    \left\{\left .\ttwo abcd \,  \right| \  ~a,\, b,\, c,\, d \in \Z \, , \
    a\,d\,-\,b\,c\,=\,1\,
    \right\}\, ,
\end{equation}
  the group of substitutions
generated by $\tau \rightarrow {\tau} + 1$ and $\tau \rightarrow
\frac{-1}{\tau}$.

\subsection{Weil's Converse Theorem (1967)}
\label{Weil}

A. Weil in 1967 completed  Hecke's theory by similarly
characterizing modular forms for  {\em congruence} subgroups, such
as
\begin{equation}\label{gammandef}
        \Gamma_0(N) \ \ = \ \
        \left\{ \left .\left( \begin{array}{cc} a & b \\
                                                  c & d
\end{array} \right) \in SL(2,\Z) \, \right| \ c \, \equiv 0 \!\!\!\!\pmod{N}
\, \right\}.
\end{equation}
(These subgroups in general have many generators, whereas Hecke's
Theorem deals with modular forms only for the groups generated by
$\tau \mapsto \tau +  \lambda$ and $\tau \mapsto -
\frac{1}{\tau}$.) Weil's breakthrough was to {\em twist} the
series $\phi(s)$ by {\em Dirichlet characters}.  Recall from
\secref{Hecke I} that a Dirichlet character modulo $r$ is a
periodic function $\chi:\Z\rightarrow \C$ which is completely
multiplicative (i.e. $\chi(nm)=\chi(n)\chi(m)$), and satisfies
$$\chi(1) \ \ = \ \ 1 \, , \ \ \  \ \ \ \ \chi(n) \ \ = \ \ 0,
\ \hbox{ if }(n,r)>1.$$  Given a Dirichlet character $\chi$ modulo
$r$ and a proper multiple $r'$ of $r$, one may form a Dirichlet
character $\chi'$ modulo $r'$ by setting $$\chi'(n) \ \ = \ \
  \begin{cases}
\,    \chi(n) & ,~(n,r')\, = \, 1, \\
 \,   0 & ,~\text{otherwise}.
  \end{cases}$$
Such a character $\chi'$ obtained this way is termed {\em
imprimitive}, and one which is not, {\em primitive}.  The
importance of primitive characters is that their functional
equations are simpler (see \cite{Daven}).  Weil's converse theorem
gives a condition for modularity under $\G_0(N)$ in terms of the
functional equations of Dirichlet series twisted by primitive
characters:

\begin{thm}\label{weilthm}(Weil \cite{Weil})
 Fix positive integers $N$ and $k$, and suppose $L(s) \, = \, \sum_{n=1}
^{\infty} \ a_{n} \,n^{-s}$ satisfies the following conditions:

\begin{itemize}
 \item{(i)} $L(s)$ is absolutely convergent for $\Re s$ sufficiently large;
\item{(ii)} for each primitive character $\chi$ of modulus $r$
with $(r,N)=1,$ $$\L(s,\chi) \ \ = \ \ (2{\pi})^{-s}\ {\Gamma(s)}\
\sum_{n=1} ^{\infty} \ a_{n}\,\chi(n)\,n^{-s}$$ continues to an
{\bf E}ntire function of $s$, {\bf B}ounded in {\bf V}ertical
strips; \item{(iii)} Each such $\L(s,\chi)$ satisfies the {\bf
F}unctional {\bf E}quation
\begin{equation}\label{weilcondiii}
 \L(s,\chi) \ \ = \ \ w_\chi \ r^{-1} \ (r^{2}N)^{\frac{k}{2}-s}
 \    \L(k-s,\bar\chi) \, ,\end{equation} where
$$w_\chi \ \ = \ \ i^k\,\chi(N)\,g(\chi)^{\,2}$$ and
    the {\em Gauss sum} $$g(\chi) \ \ = \ \
    \sum_{n\!\!\! \pmod r}\chi(n)\  e ^{\,2\,\pi\, i \,n/r}\,.$$
\end{itemize}
Then $f(z) \, = \, \sum_{n=1}^\infty \ a_n \, e^{\,2\,\pi\,  i \,
n \, z}$ belongs to the space of modular forms for $\Gamma_0(N)$
(i.e.
$$f\(\f{a\, z+b}{c\, z+d}\) \ \ = \ \ (c\, z+d)^k \ f(z) \  ,~~~ \ \
\text{for all}  \ \ \left( \begin{array}{cc} a & b \\ c & d
\end{array} \right) \, \in  \, \G_0(N) \ ,$$
and satisfies a holomorphy condition at its ``cusps'' analogous to
property (iii) in the definition of $M(\l,k,C)$ in \secref{Hecke
II}).
\end{thm}

Note that the trivial character (with $\chi(n)\equiv 1$) is
primitive, so the statement includes the $L$-functions used in
\thmref{heckecon}. Property (iii) is certainly satisfied if $L(s)$
is the $L$-function of a modular form, as can be shown using a
slight variant of Hecke's argument used in proving
\thmref{heckecon}. For the obvious reason, we refer to this
Theorem as ``Weil's converse to Hecke Theory''. For a proof, see
\cite{Bump}, \cite{Iwaniec}, or \cite{Ogg}.

\subsection{Maass Forms (1949)}\label{maass sec}

In addition to the holomorphic modular forms on the complex upper
half plane $\U$, there are the non-holomorphic modular forms
introduced by Maass \cite{maass}. These are equally important, but
far more mysterious. The literature
 has slight differences in the terminology, but for us a Maass
form will be a non-constant eigenfunction of Laplace operator
$\D=-y^2\(\f{d^2}{dx^2} + \f{d^2}{dy^2}\)$ in $L^2(\G\backslash
\U)$, where $\G$ is a discontinuous subgroup of $SL(2,\R)$, e.g. a
congruence subgroup. The laplacian condition replaces the
holomorphy condition here. In contrast to the holomorphic modular
forms, all of which have constructions and geometric
interpretations, the vast majority of Maass forms lack
constructions or identification.  Their mere existence is so
subtle that Selberg invented the trace formula \cite{sel1956}
simply to show that they exist for $\G=SL(2,\Z)$!  In fact,
deformation results such as those of Phillips-Sarnak  and Wolpert
\cite{philsar1,philsar2, philsar3,sarps,wolp1,wolp2} demonstrate
that Maass forms are scarce for the generic discrete subgroup
$\G\subset SL(2,\R)$. For this reason we shall stick to congruence
subgroups $\G$ for the rest of this exposition.

For now, consider a Maass form $\phi$ for $\G=\G_0(N)$ (for
simplicity the reader may take $N=1$ and  $\G=SL(2,\Z)$).  The
Fourier expansion of $\phi$ is given by
\begin{equation}\label{intmaass}
 \phi(x+iy) \ \ = \ \    \sum_{n\neq 0}
    a_n\ \sqrt{y}\ K_\nu( \, 2\,\pi\,|n|\,y \, )\ e^{2\, \pi\,i\,n\,x} \, ,
\end{equation}
where $a_n$ are coefficients and
 $K_\nu(t)$ is the $K$-Bessel function
\begin{equation}
\label{kdef}
  K_s(z) \ \  =
  \ \  \frac{\pi}{2} \,  \frac{{I_{-s}(z)}- I_{\, s}(z)}{\sin{\pi s}}
  \  \   \ , \ \ \ \  \ \ \ \ I_s(z) \ \ =  \ \ \sum_{m=0}^{\infty} \
  \frac{({\frac{z}{2}})^{\,s+2m}}
       {m! \  \Gamma(s \, + \, m \, + \, 1)} \, .
  \end{equation}
The parameter $\nu$ is related to the Laplace eigenvalue of $\phi$
by $\l=1/4-\nu^2$, where $\D\phi=\l\phi$. Hecke's method was
extended by Maass to obtain the analytic continuation and
functional equations of the $L$-functions
$L(s,\phi)=\sum_{n=1}^\infty a_n n^{-s}$ of Maass forms on
$\G_0(N)$ through the integral $\int_0^\infty
\phi(i\,y)\,y^{s-1/2}\,\f{dy}{y}$. When $\G=SL(2,\Z)$, for
example, this integral  is unchanged by the substitution $s\mapsto
1-s$.  Maass also proved a converse theorem for his Maass forms
for $\G=SL(2,\Z)$; see the comments at the end of the section.

\subsection{Hecke Operators}\label{heckopsec}

The Euler product structure of the Riemann $\zeta$-function has an
analog for modular form $L$-functions through Hecke operators. For
any positive integer $n$, the Hecke operator
\begin{equation}\label{tndef}
T_n(f)(z) \ \ = \ \ \f{1}{n} \ \sum_{a \,d \ = \ n }  \ a^k
\sum_{0  \, \le \,  b \,  < \,  d} \ f\(\f{a\, z+\,b}{d}\)
\end{equation}
preserves the space of modular forms of weight $k$ for $\G_0(N)$,
so long as $n$ and $N$ are relatively prime.  The same formula
applies to Maass forms when $k=0$ and the prefactor $\f 1n$ is
replaced by $\f{1}{\sqrt{n}}$. A few other operators are used as
well, to take into account symmetries of $\G_0(N)$ by which
modular forms can be ``diagonalized.''  In addition to being
eigenfunctions of a differential operator (i.e. either the
Cauchy-Riemann operator $\bar{\partial}$ or the laplacian $\D$), a
basis of modular forms or Maass forms  can be chosen among
eigenfunctions of the Hecke operators as well. As a result,
identities amongst the coefficients can be proven.  These are
nicely expressed as factorizations of the  $L$-functions of
modular forms. For example, when $\G=SL(2,\Z)$ the $L$-series of a
holomorphic form of weight $k$ factors as
\begin{equation}\label{eulerholom}
    L(s) \ \ = \ \   \sum_{n=1}^\infty \,a_n\, n^{-s}     \ \ = \ \
    \prod_{p} \(1 \, - \, a_p \, p^{-s} \, + \, p^{\,k-1-2s}\)\i\, ,
\end{equation}
a formula which remains valid for Maass forms if $k$ is taken to
be 1.

We end this section with some  remarks about the Converse
Theorem~\ref{weilthm}.  Maass observed that Hecke's argument for
\thmref{heckecon} applies to his Maass forms for
 $\G=SL(2,\Z)$ as well, but this
method does {\em not} prove a converse theorem for $\G_0(N)$ for
$N$ large. The reason for this is that the group $\G_0(N)$ can
have many generators, which are not accounted for by simply one
functional equation alone. Interestingly, Conrey and Farmer
\cite{conreyfarm} have found that by using Hecke operators, a
converse theorem can be proved for a surprisingly large range of
$N$ using only a single functional equation.  In another
direction, Booker \cite{Booker} has recently discovered that the
converse theorem requires only a single functional equation when
it is specialized to the $L$-functions coming from {\em Galois
representations}, regardless of how large $N$ is. It is an open
question whether or not Weil's argument applies to Maass forms.  A
key point for Weil is that radially symmetric holomorphic
functions are necessarily constant; this is not true in the
non-holomorphic case because there are spherical functions (formed
by radially-symmetrizing $\Im\!\!(z)^s$), and so Weil's argument
does directly apply. However, there is nevertheless an applicable
converse theorem due to Jacquet and Langlands, which we will come
to in \secref{JL}.

\section{$L$-functions from Eisenstein Series (1962-)}\label{Eisen I}

In the last section we saw the Mellin transform provided a
connection between holomorphic modular forms, and certain
Dirichlet series generalizing $\zeta(s)$. Another quite different
connection comes from a family of non-cuspidal modular forms, the
Eisenstein series
\begin{equation}\label{Gkdef}
    G_k(z) \ \ = \ \
    \sum_{(m,n)\,\in\,\Z^2-\{0,0\}}
    \f{1}{(m\,z+n)^{k}}\,, \ \ \ \ \ k\text{~even}, ~\ge 2.
\end{equation}
It is not difficult to show that $G_k(z)$ is a holomorphic modular
form of weight $k$ for $SL(2,\Z)$.  Via a Poisson summation
argument over $m$ one can obtain the Fourier expansion
\begin{equation}\label{Gkexp}
    G_k(z) \ \ = \ \    2\, \zeta(k) \ + \
    \f{2\, ( \, 2\, \pi \, i\, )^{\,k}}{(k-1)!}\
     \sum_{n=1}^\infty \,
    \sigma_{k-1}(n)\ e^{2\,\pi \,i\, n\, z}\, ,
\end{equation}
where $\sigma_{k-1}(n)$ is defined in terms of the divisors of $n$
by
$$\sigma_{t}(n) \ \ = \ \ \sum_{d|n} \, d^{\,t}$$
 (see \cite[Section
7.5.5]{Serre} for details).
 The appearance of $\zeta(k)$ here is the first example of a
very general phenomena, which  ultimately leads  to the
Langlands-Shahidi method (\secref{LaSh}).  In the next section, we
will describe the generalized {\em non-holomorphic} Eisenstein
series considered by Selberg, and their connection to the analytic
properties of the Riemann $\zeta$-function throughout the complex
plane -- not just at special integral values alone.

\subsection{Selberg's Analytic Continuation}
\label{Sel}

Selberg's method \cite{Sel} can be used to obtain the analytic
continuation and functional equations of the $L$-functions that
arise in the ``constant terms'' of Eisenstein series. We shall
sketch a form of it in the classical case of the upper half plane
$\U=\left \{z=x+iy \mid y>0 \right \}$, and the simplest-possible
Eisenstein series. Here we shall summarize the main steps
involved; the details can be found in \cite{Ku,Borel}. We will
turn to the general case in \secref{LaSh}.

Define
\begin{equation}\label{ezsdef}
    \aligned
E(z,s)\ = \ \  & \ \
 \f 1{2}~~~~
 \sum_{\stackrel{\scriptstyle (m,n)=\Z^2-\{0\}}{\scriptstyle gcd(m,n)=1}}
   \ \frac {y^s}{|m\,z+n|^{2s}} \  \\ = \  \ &  \ \f 12 \, \f{1}{\zeta(2\,s)}
\sum_{(m,n)=\Z^2-\{0\}}\frac {y^s}{|m\,z+n|^{2s}} \,
  ,
    \endaligned
\end{equation}
for $z \in \U$ and  $\sigma=\Re{s}>1.$ This series converges
absolutely and uniformly in any compact subset of the region $\Re
s >1$, and is the first example of a non-holomorphic
 Eisenstein series.  Very importantly, $E(z,s)$ is unchanged by the
 substitutions $z\mapsto \f{az+b}{cz+d}$ coming from any matrix in
(\ref{sl2zdef}).   Selberg considers  the problem of analytically
continuing $E(z,s)$ with respect to $s$ to obtain another
functional equation, as we shall now explain. (Actually, Selberg
had several different arguments to do this, but they mainly appeal
to spectral theory to obtain the important properties of analytic
continuation and functional equation of Eisenstein series.)

To motivate the statement of the functional equation, let us first
consider the Fourier expansion of $E(z,s)$. It is given by

\begin{equation}\label{eisfourexp}
E(z,s) \ =  \ E(x \, + \, i \, y \, , \, s) \  = \ \sum_{m\in\Z} \
a_m(y,s)\,e^{\,2\,\pi \, i \, m \, x}
\end{equation}
where $e(x)=e^{\,2\,{\pi}\,i\,x}$, and

$$a_m(y,s) \ = \ \int_{0}^{1} \ E(x+iy,s)\,
e^{-2\,\pi \, i \, m \, x} \, dx\, .$$ We shall need here only the
coefficients $a_0$ and $a_1$. If one computes directly for $\Re
s>1,$ using the ``Bruhat decomposition''\footnote{The Bruhat
decomposition states that all matrices $\ttwo abcd  \in SL(2,\Z)$
with $c\neq 0$ may be written as products  $\ttwo 1r{0}1
\ttwo{\a}{\b}{\g}{\d}\ttwo 1s{0}1$, where:  $r$ and $s$ range over
 $\Z$; $\g$ over $\Z-\{0\}$; $\d$ over $(\Z/\g\Z)^*$; and $\a$
and $\b$ are any two integers (which depend on $\g$ and $\d$, of
course) satisfying
 $\a\d-\b\g=1$. }  for $SL(2,\Z)$
  and recalling (\ref{kdef}),
one obtains

\begin{equation}\label{a0}
 a_0(y,s) \ \ = \ \  y^s \ + \ \phi(s)\,y^{1-s} \end{equation}
and \begin{equation}\label{a1}a_n(y,s) \ \ = \ \ 2\ \f{\sqrt{y} \
K_{s-\scriptstyle{\frac{1}{2}}}\,(2\,\pi\,|n|\,y)}{\pi^{-s}\,
\G(s)\,\zeta(2s)}\ |n|^{s-1}\ \sigma_{1\,-\,2\,s}(n) \, ,
\end{equation} with

    $$\phi(s) \ \ = \ \
     {\pi}^{1/2}\,\frac{\Gamma(s-\frac{1}{2})}{\Gamma(s)}
                \frac{\zeta(2s-1)}{\zeta(2s)}
                \ \ = \ \ \f{\xi(2s-1)}{\xi(2s)}$$
    (see
\cite{Bump} for details). In general, $\phi(s)$ is called the
``constant term'' or ``scattering'' matrix of $E(z,s)$.

Having described the Eisenstein series $E(z,s)$, we now state and
prove Selberg's theorem:

\begin{thm}\label{eisfe}(Selberg -- see \cite{Sel})
 $E(z,s)$ has a meromorphic continuation to the whole complex s-plane,
and satisfies the functional equation
\begin{equation}\label{eisesfe}
E(z , s) \ \ = \ \ \phi(s) \ E( z, 1-s)\,.
\end{equation}

\end{thm}

\vspace{.5cm}

 {\bf (A Misleading) Proof:}
\thmref{zetaebvfe} and (\ref{a0}-\ref{a1}) can then be applied to
show that each term in the Fourier expansion $$E(z,s) \ \ = \ \
\sum_{n\in \Z} \, a_n(y,s)\,e^{\, 2\, \pi \, i \, n \,  x}$$ is
meromorphic and satisfies the functional equation (\ref{eisesfe}).
The sum converges rapidly because $K_s(y)$ decays exponentially as
$y\rightarrow \infty$. Hence the whole sum is meromorphic on $\C$,
and satisfies (\ref{eisesfe}). \bx

We wrote that the above proof is ``misleading'' because, although
it demonstrates a connection to \thmref{zetaebvfe}, in practice it
has turned out to be much more fruitful to reverse the logic --
and conclude properties of $L$-functions from those of Eisenstein
series! Indeed, \thmref{eisfe} can be proven using spectral
theory, and  even in a very non-arithmetic setting (see
\cite{Borel,Cohen-Sarnak,Ku} for more details). The reader may
already have noticed a similarity between the Fourier expansion of
Eisenstein series in (\ref{eisfourexp}-\ref{a1}), and those of
Maass forms in (\ref{intmaass}).  In fact, the Eisenstein series
$E(x+iy,s)$ is an eigenfunction of the Laplace operator $\D =
-y^2\(\f{d^2}{dx^2} + \f{d^2}{dy^2}\)$, with eigenvalue $s(1-s)$.
 A main point in arguing the functional equation is {\em Maass' lemma},
  which ultimately implies
that because $E(z,s)$ and $E(z,1-s)$ share the Laplace eigenvalue
$s(1-s)$, the two must be multiples of each other.  The ratio can
be found to be $\phi(s)$ by the inspecting the constant term
$a_0(y,s)$, and so the functional equation (\ref{eisesfe}) can be
proven without knowing $\zeta$'s functional equation $\xi(s) =
\xi(1-s)$.

Of course, Selberg proved his \thmref{eisfe} in much greater
generality than we have stated. Our point is that the analytic
continuation and functional equation for the Eisenstein series
furnish an analytic continuation and functional equation for the
Riemann $\zeta$-function. To analytically continue $\zeta(s)$,
basically ``the constant term'' is enough: reading through the
spectral proof of the analytic continuation of $\phi(s)$ for
$E(z,s)$, one demonstrates that $\xi(s)$ is holomorphic
everywhere, save for simple poles at $s=0$ and 1. To get the
functional equation, we need to consider the \emph{non-trivial}
Fourier coefficient $a_1(y,s)$. \thmref{eisfe} yields
\begin{equation}\label{zfefromeis}
\aligned
 \f{2\sqrt{y}\,K_{s-1/2}(2\,\pi \,y)}{\xi(2s)}\ \ = & \ \ a_1(y,s)  & \ \
 \\= \ &
\frac{\xi(2s-1)}{\xi(2s)}  \, a_1(y,1-s) \ = \
\f{\xi(2s-1)}{\xi(2s)}
  \f{2\sqrt{y}\,K_{1/2-s}(2\,\pi \,y)}{\xi(2-2s)}\, ;
\endaligned
  \end{equation}
then, using $K_s= K_{-s}$ and setting $s=\frac{1+s'}{2}$, we have

$$\xi(s') \ \ = \ \ \xi(1-s')\, ,$$
exactly the {\bf F}unctional {\bf E}quation for $\zeta(s')$.
Incidentally, the same analysis applied to the general Fourier
coefficient $a_n(y)$ from (\ref{a1}) does not give any additional
information (this is because the extra factor
$|n|^{s-1}\sigma_{1-2s}(|n|)$ already obeys the functional
equation). {\bf B}oundedness in {\bf V}ertical strips is another
matter, which we will return to in \secref{eisbvsec}. Selberg's
work on $GL(2)$ was extended by Langlands \cite{La1,La2} to cover
Eisenstein series on general groups, where the analysis is much
more difficult.  This forms the basis of the Langlands-Shahidi
method, the topic of \secref{LaSh}.

\section{Generalizations to Adele Groups} \label{adelesec}

In the remaining sections of the paper, we will revisit the
techniques and topics of the earlier sections, but in the expanded
setting of automorphic forms on groups over the adeles.  The
adeles themselves enter as a language to keep track of the
arithmetic bookkeeping needed for complicated expressions, such as
the computations over general number fields in \secref{Hecke I}.
They are convenient even in the simplest examples when the ground
field is $\Q$.  For instance, we shall see in the next section how
Tate's thesis naturally produces the Euler product formula for the
Riemann $\zeta$-function:
$$\zeta(s) \ \ = \ \ \prod_{p}\(1-p^{-s}\)\i \ \ = \ \ \prod_p
\(1+p^{-s}+p^{-2s}+\cdots\),$$ a formula which itself is a
restatement of the unique factorization theorem for integers. They
will be useful in \secref{LaSh} for computations involving the
Eisenstein series for $SL(2,\Z)$ from \secref{Eisen I}. In
general, they are extremely valuable on general groups, where they
give clues for how to structure terms in large sums into an
``Eulerian" form.

The adeles and their notable features are perhaps better explained
later on, within the context of the arguments in which they are
used. Nevertheless we give the basic definitions before
proceeding. Given a rational number $x$, let
\begin{equation}\label{pval}
    |x|_p \ := \ p^{-\textstyle{{\scriptstyle{\hbox{ord}}}_p(x)}},
\end{equation}
where $\hbox{ord}_p(x)$ denotes the exponent of $p$ occurring in
the unique factorization of $x\in\Q$.  This $p$-adic valuation
defines a metric on $\Q$ by $d_p(x,y)=|x-y|_p$, and its completion
is $\Q_p$, the field of $p$-adic numbers.  More concretely, $\Q_p$
may be viewed as the formal Laurent series in $p$
\begin{equation}\label{padicform}
   x \ \ = \ \  c_k \, p^{\,k} \, + \,
   c_{k+1} \, p^{\,k+1} \,+ \, \cdots \ , \ \ \ \ \  c_k \neq 0\, , \ \
   \ \ \ k \ \in \ \Z
\end{equation}
with integral coefficients $0 \le c_j < p$;  alternatively it may
be thought of as consisting of base-$p$ expansions with only
finitely many digits to the right of the ``decimal'' point, but
perhaps infinitely many to the left. Within $\Q_p$ lies its ring
of integers, $\Z_p$, which is the completion of $\Z$ under
$|\cdot|_p$.  It may instead be viewed as the elements of $\Q_p$
as in (\ref{padicform}) which have $k \ge 0$, or those with no
digits to the right of the decimal point in their base-$p$
expansion. The $p$-adic valuation of course extends to $\Q_p$: the
absolute value of $x$ given in (\ref{padicform}) is $p^{-k}$, and
 $\Z_p = \{x\in\Q_p \, \mid \, |x|_p \le 1\}$. Similarly, the multiplicative subgroups
 are $\Q_p^*
= \Q_p\,-\,\{0\}$ and $\Z_p^* =\{x\in\Q_p \, \mid \, |x|_p = 1\}$.

The adeles are formed by piecing together all $\Q_p$ along with
$\R$, which may be viewed as $\Q_\infty$, the completion of $\Q$
under the usual archimedean absolute value. Concretely, the adeles
$\A$ are the restricted direct product of the $\Q_p$ with respect
to the $\Z_p$; that means the adeles are infinite-tuples of the
form
\begin{equation}\label{adeledef}
    a \ \ = \ \ (a_\infty;\,a_2,\,a_3,\,a_5,\,a_7,\,a_{11},\,\ldots)\,, \ \
    \ a_p \, \in \, \Q_p  \ \ \text{for all} \ \  p\, \le \, \infty
\end{equation}
such that all but finitely many $a_p$ lie in $\Z_p$.  Similarly
the ideles $\A^*$ are the restricted direct product of all
$\Q_p^*$ with respect to $\Z_p^*$.  Addition and multiplication
are defined componentwise in $\A$ and $\A^*$.  The rational
numbers embed diagonally into the ring $\A$ and play a fundamental
role, which will become apparent shortly when it appears in Tate's
thesis. The adeles, or more properly the ideles, themselves have
an absolute value; its value on $a$ in (\ref{adeledef}) is
$$|a|_\A \ \ = \ \ \prod_{p\le \infty} |a_p|_{p}.$$  Note that this
is actually a finite product, because almost all $a_p$ have
absolute value equal to one, a theme which underlies many adelic
concepts.  The diagonally-embedded $\Q^*$ consists of the ideles
with $|a|_\A = 1$.

The above construction can be generalized to an arbitrary number
field -- or even ``global field'' -- $F$ to obtain its adele ring
$\A_F$  (see \cite{Lant,ramval}). Most constructions involving
$\A_\Q$ generalize to $\A_F$, though we will mainly focus on
$F=\Q$ for expositional ease.  Adeles are usually viewed much more
algebraically and with much greater emphasis on their topology
(which we have hardly touched); our intention here is rather to
give enough background to illuminate their effectiveness in
analysis.

\section{Tate's Thesis (1950)}
\label{Tate}

In his celebrated 1950 Ph.D. thesis \cite{Tate}, J. Tate
reinterpreted the methods of Riemann and Hecke in terms of
harmonic analysis on the ideles $\A^*$ of a number field $F$.
Tate's method succeeded in precisely isolating and identifying the
contribution to the functional equation from each of the ramified
prime ideals $\mathfrak P$ not treated in the product
(\ref{heckeslseries}), a delicate problem which appeared
complicated from the perspective of Hecke's classical method.  At
the same time, Tate's method is powerful enough to uniformly
reprove the analytic continuation and functional equations of
Hecke's $L$-functions.  For this local precision, uniformity, and
flexibility, Tate's method has influenced the many adelic methods
at the forefront today. In this section we explain Tate's
construction and the role of the devices he employs, via a
comparison with Riemann's argument in \secref{Rie}.

Let us recall Riemann's integral from \secref{Rie}, after a
harmless change of variables:
\begin{equation}\label{riemannsintegral}
    \xi(s) \ = \ \pi^{-s/2}\,\G({\textstyle \f s2})\,\zeta(s) \ =
    \ \, \int_{0}^\infty x^s\,\sum_{n\neq 0} e^{-\pi \,n^2\,
    x^2}\,d^*x \, .
\end{equation}
Tate instead considers the sum over $\Z-\{0\}$ as an integral over
a disconnected group.  In order to keep the flexibility of
treating more general sums, he instead essentially integrates the
characteristic function of $\Z$ over a much larger set in his
generalized {\em $\zeta$-integral}
\begin{equation}\label{tatesintegral}
      \zeta(f,c) \ \ = \ \ \int_{{\A_\Q}^*} f(a)\,c(a) \, d^* {a} \, .
\end{equation}
Here $c(a)$ is any quasi-character of $\A^*$ -- that is, a
continuous homomorphism from $\A^*$ to $\C^*$ -- which is trivial
on $\Q^*$ (for example, we saw before that $|a|_\A$ is trivial on
$\Q^*$); $d^*a$ is the multiplicative Haar measure on $\A^*$
pieced together as a product of the local Haar measures
$d^*x_{\infty}=\f{dx}{|x|}$ and $d^*x_p$. The latter is normalized
so that $\Z_p^*$ has measure 1.  Finally, the function $f$ is
taken to be a product
\begin{equation}\label{tatefdef}
    f(a_\infty;\,a_2, \, a_3, \, a_5,\ldots) \ \ = \ \ \prod_{p\le
    \infty} f_p(a_p)
\end{equation}
of functions $f_p$ on $\Q_p$, which may depend on the
quasi-character $c$. In the simplest possibility, which is that
$c(a)=|a|_\A^s$, let us choose
$$f_p(x) \ \ = \ \ \chi_{\Z_p}(x) \
\ =  \ \
  \begin{cases}
    \,1 &,~ |x|_p \leq 1\,, \\
    \,0 &,~ \text{otherwise}\,,
  \end{cases}$$
and $f_{\infty}(x)=e^{-{\pi}\,x^2}$; then the integral
$\zeta(f,|\cdot|_\A^s)$ actually recovers Riemann's integral. This
can be seen as follows: first we may ``fold'' the integral to one
over $\Q^*\backslash \A^*$:
\begin{equation}\label{tatefold}
    \int_{\A^*} f(a)\, |a|^s_{\A}\,d^*a \ \ = \ \
     \int_{\Q^* \backslash \A^*}|a|_\A^s
     \,\(  \sum_{q\in\Q^*} f(qa) \) \, d^*a \, .
\end{equation}
The {\em strong approximation principle} states that
$(0,\infty)\times \widehat{\Z}^*$ is a fundamental domain for
$\Q^*\backslash \A^*$, where $ \widehat{\Z}^* = \prod_{p<\infty}
\Z_p^*$.  It is easy to see that $f(qa)\equiv 0$ on this
fundamental domain unless the rational $q$ is actually an integer,
for otherwise, the $p$-adic valuation $|qa|_p =|q|_p > 1$ for any
prime $p$ in the denominator of $q$.  Thus the role of the $f_p$
is to select the integers amongst $\Q$, and (\ref{tatefold})
becomes
\begin{equation}\label{tatefold2}
 \int_{(0,\infty)\times \widehat{\Z}^*}|a|_\A^s
     \,\(  \sum_{n\neq 0} f(na) \) \, d^*a \, .
\end{equation}
Now $f_p((na)_p)\equiv 1$ for all $p<\infty$, and so the integrand
is independent of the $\widehat{\Z}^*$ factor, which has volume 1
under the Haar measure. Now (\ref{tatefold2}) amounts to
\begin{equation}\label{tatefold3}
    \int_0^\infty |a_\infty|^{\,s} \,\sum_{n\neq 0}
    e^{-\pi \, n^2\,a_\infty^2} \,\ d^*a_{\infty} \, ,
\end{equation}
i.e. (\ref{riemannsintegral}).  Thus Tate's and Riemann's
integrals match for $\zeta(s)$.

At the same time, the global integral on the lefthand side of
(\ref{tatefold}) factors as a product
\begin{equation}\label{tateprod}
    \prod_{p\le \infty} \, \int_{\Q_p^*} f_p(x)\ |x|_p^s\ d^*x_p \ \
    =
    \ \ \left(\int_{\R} e^{-{\pi}|x|^2} |x|^{s} \ \f{dx}{|x|}
\right) \, \cdot \, \prod_{p}\, \int_{\Z_p}
   { |x_p|_p^{s}} \  d^*x_p \, .
\end{equation}
The integral over $\R$ gives $\pi^{-s/2}\G(\f s2)$, and the
$p$-adic integral may actually be broken up over the ``shells''
$p^k\,\Z_p^*=\{|x_p|_p = p^{-k}\}$, $k\ge 0$, to give the
geometric series $\sum_{k=0}^\infty p^{-k\,s} = (1-p^{-s})\i$.
This gives the Euler product formula for $\zeta(s)$, along with
its natural companion factor $\pi^{-s/2}\G(\f s2)$ for $p=\infty$
-- in other words, the completed Riemann $\xi$-function.

We should note that the role of the adelic absolute value (and in
particular that its value is 1 on $\Q$) corresponds to the change
of variables $x\mapsto x/n$ in the classical picture. In general
for a global field $F$, we may  write the quasi-character $c(a)$
in the form $c_0(a)|a|^s$, where $c_0: {\A^*} \rightarrow \C^*$ is
a character of modulus 1. Then $c_0(a)$ corresponds to $\chi$, a
``Hecke character'' for $F$ (\secref{Hecke I}), and $\zeta(f,c)$
differs from

      $$L_F(s,\chi) \ \ = \ \  \prod_{{\mathfrak P}}
       \( 1\,-\,{\chi}({\mathfrak P})\,(N{\mathfrak P})^{-s}\)^{-1}$$
(where ${\mathfrak P}$ now runs over {\em all} prime ideals of
$F$) by only a finite number of factors.

In this idelic setting, Tate uses a  Fourier theory and Poisson
summation formula on the ring of adeles $\A$, and proves the
elegant functional equation
\begin{equation}\label{tatelegent}
        \zeta(\,  f\, , \,  c\,  ) \ \ = \ \ \zeta( \,  \widehat{f} \, ,\, \widehat{c}\,)\, ,
\end{equation}
where $\widehat{f}$ is the ``adelic Fourier transform'' of $f$ and
${\widehat{c}}(a)={\overline{c_0(a)}}{|a|^{1-s}}$. The functional
equation for $L_F(s,\chi)$ may be extracted from this. To
illustrate with our example of the Riemann $\zeta$-function,
recall that we had taken $c_0$ to be identically equal to 1, and
in fact our $f=\widehat{f}$, so that
   $$\zeta(f_0,|\cdot|^s)  \ \ = \ \
   \pi^{-\frac{s}{2}}\, \Gamma({\textstyle \frac{s}{2}})\,
   \prod_{p}\,(1-p^{-s})\i \ \ = \ \ \xi(s)\, .$$
The functional equation {\bf FE} is then immediate from
(\ref{tatelegent}). Tate's method of course also yields the {\bf
E}ntirety and {\bf B}oundedness in {\bf V}ertical strips.

\section{Automorphic forms on $GL(n)$}\label{glnautsec}

Thus far we have seen two types of $L$-functions: The Riemann
$\zeta$-function and its cousins that are treated in the
Riemann-Hecke-Tate theory (\secref{Rie}), and the $L$-functions of
modular forms in Hecke's (other) theory (\secref{Hecke II}).  We
now understand these $L$-functions to be part of a family, the
$L$-functions of automorphic forms on $GL(n,\A)$.  The integrals
in Tate's thesis are over $\A^*$, which is just $GL(1,\A)$, and
the quasi-characters are viewed as automorphic forms on
$GL(1,\Q)\backslash GL(1,\A)$.
 We shall now explain how to view the modular
forms we saw in \secref{Hecke II} as automorphic forms on
$GL(2,\Q)\backslash GL(2,\A)$.  This leads to two major
generalizations: first to a general number field (or indeed even a
global field) $F$ instead of $\Q$, and second to an arbitrary
reductive algebraic group $G$ instead of $GL(1)$ or $GL(2)$.

To recap from \secref{Hecke II}, a holomorphic modular form of
weight $k$ for $\G=SL(2,\Z)\subset SL(2,\R)$ is a holomorphic
function on the complex upper half plane $\U$ such that
\begin{itemize}
\item $f\(\f{az+b}{cz+d}\)\  =   \ (cz+d)^k\,f(z)$ for all $\ttwo
abcd \in \G$ \item $f(z)$ has a Fourier expansion $f(z)=\sum_{n\ge
0} c_n\,e^{2\, \pi \, i \, n \, z}$.  In addition $f$ is a cusp
form if $c_{\,0}=0$.
\end{itemize}
The above definition of course extends to more general groups
$\G$, such as the congruence subgroups in (\ref{gammandef}).
Before considering $f$ as a function on $GL(2,\A)$, we must first
explain how to consider $f$ as a function on $GL(2,\R)$, or even
$SL(2,\R)$.  Indeed, there is a correspondence between holomorphic
modular forms $f$ of weight $k$ for $\G\backslash \U$, and certain
functions $F$ on  $\G\backslash SL(2,\R)$ defined via the
following relations:
\begin{equation}\label{htosl2relation}
\gathered
  F\ttwo abcd \ \ = \ \
    f\(\f{a\,i+b}{c\,i+d}\)(c\,i+d)^{-k}\, ,
\\
f(x+i\,y) \ \ = \  \ y^{-k/2}\,F\(\ttwo 1x{}1
\ttwo{\sqrt{y}}{}{}{1/\sqrt{y}}\) \, .
\endgathered
\end{equation}
(We leave matrix entries blank if they are zero.) The key reason
for this correspondence is that $\U$ is isomorphic to the quotient
$SL(2,\R)/SO(2,\R)$. For more details and a precise
characterization of $F$, see \cite{Bump} or \cite{gelborange}.

Aside from the holomorphic modular forms, the most significant
automorphic forms on $\U$ are the non-holomorphic Maass forms:
non-constant, $L^2$ Laplace eigenfunctions on the quotient
$\G\backslash \U$. We described these in \secref{maass sec}.
Because of the identification $\U$ $\cong$ $SL(2,\R)/SO(2,\R)$,
Maass forms can directly be viewed as functions on $\G\backslash
SL(2,\R)$.

Now that we view the holomorphic and Maass modular forms on the
group $G=SL(2,\R)$, wide generalizations are possible, and
techniques from representation theory may be applied.  The group
$G$ acts on $L^2(\G\backslash G)$ by the right regular
representation, which is translation on the right:
\begin{equation}\label{rrr}
    [\rho(g)f](h) \ \ = \ \ f(hg) \, .
\end{equation}
The study of automorphic forms on $\G\backslash G$ now becomes
understanding the decomposition of the very large (and highly
reducible) representation $\rho$ into irreducible components. This
is the starting point for the notion of ``automorphic
representation," but for that we first need to delve more into the
arithmetic nature of $\G$, and consider $G$ adelically.

In addition to the action on the right, left-translation by
rational matrices is very important in many constructions in
automorphic forms.  We have, therefore, also the left regular
representation:
\begin{equation}\label{lrr}
    [\l(g)f](h) \ \ = \ \ f(g\i h)\, ,
\end{equation}
which maps $L^2(\G\backslash G)$ to $L^2 (g\G g\i \backslash G)$.
In general this moves automorphic forms for one congruence
subgroup $\G$ to those on a conjugate, which may be wildly
different.  For this reason it is natural to act on the left only
by rational matrices $g$, so that the conjugate of $\G$ is still
closely related to a congruence subgroup.\footnote{See
\cite{margulisbook} for a thorough explanation.} In fact, many
fundamental constructions (such as Hecke operators) require action
by rational matrices $g$ which lie in $GL(2,\Q)$, but not
$SL(2,\Q)$.  It is for this reason that we will consider adelic
automorphic forms on $GL(2,\A)$, not $SL(2,\A)$, though a theory
exists for that group as well. Because $GL(2,\R)$ is one dimension
larger than $SL(2,\R)$,  we technically need to consider
$L^2_\omega(Z \G\backslash GL(2,\R))$ where $Z$ is the center of
$GL(2,\R)$ (=scalar multiples of the identity matrix).  Here
$\omega$ is a central character (that is, a character of $Z$) and
this $L^2$ space consists of functions on $G$ which transform by
$Z$ according to $\omega$, but which are otherwise
square-integrable on the quotient $Z \G\backslash GL(2,\R)$.  As a
minor technicality, we will now consider $\G=GL(2,\Z)$ instead of
$SL(2,\Z)$ to make the picture more uniform.  The setup of this
paragraph works equally well for $\G=GL(n,\Z)$ and $G=GL(n,\R)$.

Finally we now come to the adeles.  The adelic group $GL(n,\A)$ is
the product of $GL(n,\R)$ with $GL(n,\A_f)$, the direct product of
all $GL(n,\Q_p)$ with respect to their integral subgroups
$GL(n,\Z_p)$.  We have already seen -- at least in the case $n=2$
-- that the first factor, $GL(n,\R)$, acts on automorphic
functions on $Z\G\backslash G$ on the right, and that rational
matrices act on the left.  Just as with $\A^*=GL(1,\A)$ in Tate's
thesis, there is a version of the strong approximation theorem for
$GL(n,\A)$. It states that $GL(n,\A) = GL(n,\Q) GL(n,\R) K_f$,
where $K_f = \prod_{p<\infty} GL(n,\Z_p)$.   We now define an
action of $GL(n,\A_f)$ on the left that extends the action of
$GL(n,\Q)$:
\begin{equation}\label{adelicleftaf}
    [\l(g_f)F](h)  \ \ = \ \
    [\l(\g)F](h) \ \ = \ \  F(\g\i
    h) \, ,
\end{equation}
 where $\g\in GL(n,\Q)$ is the factor guaranteed by the
  strong approximation theorem in writing $g_f \in GL(n,\A_f)\subset
  GL(n,\A)$ as a product.
This definition is well defined, because any two close
``approximants'' $\g$ must be related by a multiple of an integral
matrix, and $F$ is presumed to be invariant under $GL(n,\Z)$.
Roughly speaking, the topology on the adeles is given in terms of
a basis of products of $GL(n,\Z_p)$ and finite index subgroups,
which are related to congruence groups.  Thus  the adelic topology
aligns with the invariance of $F$ under matrices in a congruence
subgroup.

Unifying these actions leads to the notion of adelic
representations and adelized automorphic forms, where $F$, instead
of being a function on $GL(n,\R)$ alone, is padded with extra
variables. Namely, we have an adelic function
\begin{equation}\label{paddedF}
    F_\A(g_\infty\,;\, g_2\, ,g_3\, ,g_5\,\ldots)
\end{equation}
such that almost all $g_p$ lie in $GL(n,\Z_p)$, and
\begin{equation}\label{FAdef}
       F_\A(g_\infty \, ; \, g_2\, ,g_3\, ,g_5\,\ldots) \ \ = \ \
       [\l(g_2)\,\l(g_3)\,\l(g_5)\,\cdots\l(g_p)\,\cdots
       F](g_\infty) \, ,
\end{equation}
where the number of $\l(g_p)$'s that act in an nontrivial way is
finite, and their actions for various $p$'s commute with each
other. The adelized function $F$ has the properties
\begin{equation}\label{Fprops}
   F(g_\infty) \ \ = \ \ F_\A(g_\infty \, ; \, 1,1,1,\cdots)
    \ \ \ \ \ \
   \text{and}\ \  \ \  F_\A(\g\,g) \ \ = \ \ F_\A(g) \, ,
\end{equation}
for any diagonally embedded rational matrix $\g$.  The center $Z$
and central character $\omega$ have  analogous adelic versions,
which are related to automorphic forms on $GL(1,\A)$, in fact.
 The right
regular representation now acts on $GL(n,\A)$ by the formula
(\ref{rrr}), but note that this right action of the factor
$GL(n,\A_f)$ is really a left action on $GL(n,\R)$.

 This leads us
to our final version of automorphic representation: an irreducible
subrepresentation of the action of the right regular
representation $\rho$ on $GL(n,\A)$ on $L^2_\omega(Z_\A \,
GL(n,\Q)\backslash GL(n,\A))$.  The constituents of these
subspaces are generalizations of the automorphic forms we
encountered previously. All have classical counterparts as
functions with  transformation properties for various congruence
subgroups $\G$ of $GL(n,\Z)$, but the adelic version provides a
uniform framework. The role of $\G$ itself is replaced by
right-invariance under finite index subgroups $K_f'$ of $K_f =
\prod_{p<\infty} GL(n,\Z_p)$; here $\G = \{\g_\R \mid \g\in
GL(n,\Q),\g_f \in K_f'\}$, where $\g_\R$ and $\g_f$ denote the
projections of $\g\in GL(n,\A)$ to the factors $GL(n,\R)$ and
$GL(n,\A_f)$, respectively. We note that forms for
 various conjugate subgroups all fall into the same
irreducible adelic automorphic representation, as do ``newforms.''
More importantly, the theory of Hecke operators
(\secref{heckopsec}) for powers of a prime $p$ (e.g. the $T_{p^k}$
from \secref{heckopsec}) can be recast as the study of the action
of $GL(n,\Q_p)$, which allows the powerful representation theory
of this $p$-adic group to be used. In general, the adelic
framework unifies many constructions in automorphic forms, and
explains their effectiveness.  Better yet, it provides insight for
new constructions which would seem very difficult to uncover using
only the classical perspective.

\subsection{Jacquet-Langlands (1970)}
\label{JL}

In 1970, a remarkable book was published: {\it ``Automorphic Forms
on $GL(2)$''}, by H. Jacquet and R. Langlands \cite{JL}.
 The irreducible unitary
representations $\pi$ of $GL(n,\A)$ discussed above factor into
restricted tensor products $\pi \cong \otimes_{p\le \infty}
\pi_p$, where $\pi_p$ is a ``local representation'' of
$GL(n,\Q_p)$.  One can treat the case of a number field, or even
an arbitrary global field in a similar way.  For $n=2$, Jacquet
and Langlands rephrase Hecke's theory from \secref{Hecke II} using
adelic machinery, much in the way Tate reworked Riemann and
Hecke's classical arguments. In particular, they attach a global
$L$-function $\L(s,\pi)$ (a Dirichlet series times a product of
gamma factors, such as the $\pi^{-s/2}\G(s/2)$ that differentiates
$\xi(s)$ from $\zeta(s)$) to each automorphic representation of
$GL(2)$.  They prove that $\L(s,\pi)$ is ``nice,'' meaning that it
satisfies the standard properties of {\bf E}ntirety, {\bf
B}oundedness in {\bf V}ertical strips, and {\bf F}unctional {\bf
E}quation that Hecke's method yields. Secondly they give a
criteria for  any ``nice'' $L$-function of this type to come from
an automorphic representation; that is, a converse theorem.

Although the methods of group representations are new, the
underlying technique of Jacquet-Langlands is fundamentally Hecke's
method, as we shall briefly describe.  However, neither the
statement nor proof of their converse theorem is really Weil's
\thmref{weilthm}.  For example, let us return to the discussion
concluding \secref{modsec}.  Weil's proof of his converse theorem
demonstrates that only a finite number of Dirichlet characters are
required in his twisting hypothesis (ii). In fact,
Piatetski-Shapiro \cite{psbolnoi},  carefully examining this
point, discovered an important simplifying feature in the early
1970s which has become one of the most important technical devices
in today's applications.  He found that Jacquet-Langlands' proof
also requires only a finite number of twists by characters -- but
a completely disjoint set of characters from the ones Weil needed!
For a classical treatment, see \cite{razar}.

 Recall how in
\secref{Hecke II} we considered the $L$-functions of modular forms
for $SL(2,\Z)$.  The first example of a modular form whose
$L$-function is entire is {\em Ramanujan's $\D$ form}
\begin{equation}\label{ramformdef}
  \D(z) \ \ = \ \ e^{2\,\pi \,i \, z} \ \prod_{n\ge 1}
  \,(1-e^{2\, \pi \,
  i \, n \, z})^{\,24}\,,
\end{equation}
which has weight $k=12$. (See \cite{Serre} for a beautiful
exposition of $\D$ in the context of Hecke theory.) Expand the
product as $\D(z)=\sum_{n\ge 1}\tau(n)e^{\,2\,\pi\, i \,n\,z}$ and
normalize the coefficients by setting $a_n=\f{\tau(n)}{n^{11/2}}$;
in this normalization, the Ramanujan conjecture (established by
Deligne \cite{Deligne}) can be stated uniformly as $|a_p|\le 2$
for all primes $p$. Ramanujan also conjectured\footnote{though not
in the language of $L$-functions.} that the ``standard'' (i.e.
Hecke) $L$-series associated to $\D$ has an {\em Euler product}
over primes, much like $\zeta$:
\begin{equation}\label{lsdelta}
  L(s,\D) \ \ = \ \ \sum_{n\ge 1}a_n\,
  n^{-s} \ \ = \ \ \prod_{p}\,(1\,-\,a_p\,p^{-s}\,+\,p^{-2s})\i.
\end{equation}
This was proven by Mordell \cite{mordell}, and nowadays we
understand the factorization as being equivalent to the assertion
that $\D$ is an eigenfunction of the {\it Hecke operators}
(\ref{tndef}) from \secref{heckopsec} -- in particular, this is
(\ref{eulerholom}).

Let us now explain the connection between the arguments of
Jacquet-Langlands and of Hecke.  Our starting point is the Fourier
expansion of a modular or Maass  form $\phi(x+iy)$ in the variable
$x$, in which it is periodic (with period 1 in the case of
$SL(2,\Z)$, as we shall now consider). Recall that if $\phi$ is a
holomorphic cusp form of weight $k$,
\begin{equation}\label{cnholom}
\phi(z) \ \ = \ \ \sum_{n=1}^\infty \, c_n\,e^{\,2\, \pi \, i \, n
\,z}. \end{equation}
 Similarly for a Maass form we have the
Fourier expansion (\ref{intmaass}). Up to constants (and a factor
of $y^{k/2}$ in the holomorphic case, like in
(\ref{htosl2relation})), we may write these expansions as
\begin{equation}\label{firstwhit}
   \sum_{n\neq 0} \f{a_n}{|n|^{1/2}} \,W(2\,\pi \, n\, y)\,
    e^{2\,\pi\,i\,n\,x} \, ,
\end{equation}
where as before $a_n=\f{c_n}{n^{(k-1)/2}}$ in the holomorphic
case, and
\begin{equation}\label{firstwhitdef}
    W(y) \ \ = \ \ \left\{%
\begin{array}{ll}
   y^{k/2}\,e^{-y} \, , & \phi \hbox{~holomorphic,} \\
    \sqrt{|y|}\,K_\nu(|y|) \, , & \phi \hbox{~a Maass form.} \\
\end{array}%
\right.
\end{equation}
In \secref{glnautsec} we saw how both holomorphic and Maass forms
can be viewed as functions on $GL(2,\R)$.  With this point of view
we can write the corresponding function, up to constants, as
\begin{equation}\label{FfromWexpn}
F(g) \ \ = \ \ \sum_{n\neq 0} \,\f{a_n}{|n|^{1/2}}\,W\(\ttwo
n{}{}1 g\)\, , \end{equation}
 where
\begin{equation}\label{iwawhit}
W\(\ttwo 1x{}1 \ttwo {y_1}{}{}{y_2} k\) \ \ = \ \
e^{2\,\pi\,i\,n\,x}\,W(2\,\pi\,y_1/y_2) \, . \end{equation} (Here
the matrix $k$ on the lefthand side is orthogonal; all matrices in
$GL(2)$ can be written  in that form according to the {\em Iwasawa
decomposition}.) This $W(g)$ is called a ``Whittaker'' function in
connection with the special functions it is related to. It
satisfies a transformation law on the left:
\begin{equation}\label{whitleft}
    W\(\ttwo 1u{}1 g\) \ \ = \ \ e^{2\,\pi\,i\,u}\,W(g)\, ,
\end{equation}
and thus can be obtained from the integral
\begin{equation}\label{whitint}
    W(g) \ \ = \ \ \int_{0}^1 F\( \ttwo 1u{}1 g \)\,
    e^{-2\,\pi\,i\,u} \, du\,.
\end{equation}

The adelic method of Jacquet and Langlands involves incorporating
the coefficient $a_n$ into a cognate Whittaker function which
generalizes the properties (\ref{whitleft}) and (\ref{whitint}).
Consider now the adelized version $F_\A$ of $F$ defined in
(\ref{paddedF}-\ref{Fprops}), and define its adelic Whittaker
function
\begin{equation}\label{adelwhitdef}
    W_\A(g_\A) \ \ = \ \ \int_{\Q\backslash \A} F_\A\(\ttwo 1u{}1
    g_\A\)\,\psi(-u)\, du \, ,\ \ \ g_\A \,\in\,GL(2,\A) \, ,
\end{equation} where $\psi$ is a non-trivial character of $\A$
that is trivial on the subgroup $\Q$ (which we recall is
diagonally embedded into $\A$). The measure $du$ is normalized to
give $\Q\backslash \A$ measure 1.  The definition depends on the
precise choice of character, but all non-trivial characters can be
written as $\psi(qu)$ for some $q\in\Q^*$, and this $q$ can be
absorbed into $g_\A$ via the matrix $\ttwo q{}{}1$; changing
variables does not affect the measure since the adelic absolute
value $|q|_\A=1$. The result is that  $F_\A$ can be reconstructed
from $W_\A$ via the succinct formula
\begin{equation}\label{FfromW}
    F_\A(g_\A) \ \ = \ \ \sum_{q\, \in \, \Q^*} W_\A\(\ttwo q{}{}1
    g_\A\)\, ,
\end{equation}
much like (\ref{FfromWexpn}).

Now we shall make a tacit assumption that our original modular
form is a Hecke eigenform (see \secref{heckopsec}). Our Whittaker
function here, like many adelic functions, can be expressed as a
product of {\em local} Whittaker functions $W_p$ on $GL(2,\Q_p)$:
\begin{equation}\label{localwhit}
    W_\A( g_\A) \ \ = \ \ \prod_{p\le \infty} W_p(g_p)\, ,
    \ \ \ g_\A \, = \, (g_\infty\,;\, g_2\, , g_3\, , g_5\, , g_7\,
    \ldots)\,,
\end{equation}
where each $W_p$ obeys a transformation law similar to
(\ref{whitleft}).  In fact, just as with the Iwasawa decomposition
in (\ref{iwawhit}), the local Whittaker functions depend only on
diagonal matrices, and actually their value there is related to
the original Fourier coefficient by $W\ttwo {p^k}{}{}1 = a_{p^k}$.
This last fact underlies the connection between (\ref{firstwhit})
and (\ref{FfromW}): the extra adelic variables encode the  value
of the Fourier coefficients; these are very often zero, notably
when $k<0$ and the subscript is no longer an integer. This is why
the sum over $\Q$, which appears to be much larger, actually
corresponds to the sum over $\Z$.

Jacquet and Langlands use this theory beautifully to write the
global $L$-function as
\begin{equation}\label{jlintegral}
    \L(s) \ \ = \ \ \int_{\Q^*\backslash \A^*} F_\A \ttwo a{}{}1 \, |a|_\A^{s-1/2}
    \, d^*a\, .
\end{equation}
This has a functional equation $s\mapsto 1-s$, owing to the
invariance of $W$ under $\ttwo{}{-1}{1}{}$, just as in Hecke's
argument. Now by substituting (\ref{FfromW}) and collapsing the
common $\Q^*$ from the quotient and sum together (``unfolding''),
the integral
\begin{equation}\label{jlintegral2}
    \L(s) \ =  \ \int_{\A^*} W_\A\ttwo a{}{}1\, |a|_\A^{s-1/2}\,
    d^*a \ =  \ \prod_{p\le \infty} \int_{\Q_p^*}
    W_p\ttwo{a_p}{}{}1 \,  |a_p|_p^{s-1/2}\,d^*a_p
\end{equation}
splits as a product of local integrals.  The ones for $p<\infty$
separately give the local factors of an Euler product which
represents the Dirichlet series  for Hecke's $L$-function $L(s)$,
and the integral for $p=\infty$ gives the corresponding
$\G$-functions $L(s)$ must be multiplied by in order to have a
clean functional equation.  This is entirely analogous to the
situation  in Tate's thesis after (\ref{tateprod}), and the
computations are deep down identically those needed for the
classical treatment in \secref{modsec}. Details can be found in
\cite{JL,Godemondsnotes,Bump}.

We should emphasize that the method is far more general and has
strong advantages in its local precision, in that it gives a very
satisfactory treatment of the contribution to the functional
equation by each prime.  Also the technique works for congruence
subgroups, as well as over general global fields. Just as Tate's
thesis understood Riemann's $\zeta$-function in terms of $\A^* =
GL(1,\A)$, Jacquet-Langlands subsumed the theory of modular forms
and their $L$-functions through $GL(2,\A)$.
 Subsequently, efforts were underway to provide a similar theory
for general groups, most notably $GL(n,\A)$.

\subsection{Godement-Jacquet (1972)}
\label{JaGo}

Tate (see \secref{Tate}) redid Hecke (\secref{Hecke I}) by using
adeles, developing a Poisson summation formula, and working with

  $${\zeta(f,c)} \ \ = \ \ \int_{\A^*} \ f(a)\,c(a) \, d{^*}a \, .$$
R. Godement and H. Jacquet \cite{God-Jac} generalized Tate by
working with $GL(n)$ for arbitrary $n$ instead of $GL(1)$. In
particular, they proved that the global, completed $L$-functions
of automorphic
 forms on $GL(n)$ satisfy properties {\bf E, BV} and
 {\bf FE} of \secref{Intro} (actually their integral representation for $GL(2)$
 is completely different than Jacquet-Langlands and Hecke).

 We shall not pursue this avenue here,
  but will briefly describe the
 $L$-functions of cusp forms on $GL(n)$.  Recall how after
 (\ref{ramformdef}) we renormalized the coefficients
 of Ramanujan's $\D$-form by a factor of $n^{(k-1)/2}$.  This can
 be carried out for any holomorphic cusp form $f$ of weight $k$ for
 $SL(2,\Z)$, resulting in an Euler product of the same form as
 (\ref{lsdelta}).
 Writing $a_p=\a_p+\a_p\i$,  the
Euler product factors further as
\begin{equation}\label{firstlsfeul}
L(s,f)\ \ = \ \ \prod_{p} \, (1-\a_p \,p^{-s})\i \, (1-\a_p\i
p^{-s})\i.
\end{equation}  The preceding
expression is called a {\em degree two} Euler product because of
its two factors, in comparison with the degree one Euler product
$\zeta(s)\ = \ \prod_p(1-p^{-s})\i$.

The $L$-functions of cusp forms $\phi$ for $GL(n,\A_\Q)$ are Euler
products of degree $n$,
\begin{equation}\label{degreen}
    L(s,\phi) \ \ = \ \ \prod_{p}\prod_{j=1}^n \,(1 \, - \,
    \a_{p,j}\,p^{-s})\i.
\end{equation}
To form their global, completed $L$-functions, they must be
multiplied by a product of $n$ $\G$-factors,
\begin{equation}\label{glngam}
L_\infty(s,\phi) \ \ = \ \ \prod_{j=1}^n \, \G_\R(s \,+ \,\mu_j)\,
,
\end{equation}
where the $\mu_j$ are special complex parameters related to $\pi$
(for example, the $\nu$ from Maass forms in (\ref{intmaass})), and
$\G_\R(s)=\pi^{-s/2}\G(s/2)$ is again the factor which
distinguishes $\xi(s)$ from $\zeta(s)$.  The completed
$L$-function,
\begin{equation}\label{completedgln}
    \L(s,\phi) \ \ = \ \ L_\infty (s,\phi) \ L(s,\phi)
\end{equation}
is {\bf E}ntire (unless $n=1$ and $\L(s,\phi)=\zeta(s)$), {\bf
B}ounded in {\bf V}ertical strips, and satisfies the {\bf
F}unctional {\bf E}quation
\begin{equation}\label{completedglnfe}
    \L(s,\phi) \ \ = \ \ w \ Q^{1/2\, - \, s} \
    \L(1-s,\tilde{\phi})\, .
\end{equation}
Here $w$ is a complex number of modulus one, the ``conductor'' $Q$
is a positive integer (related to the congruence subgroup $\phi$
comes from), and $\tilde{\phi}$ is the ``contragredient''
automorphic form to $\phi$ (coming from the automorphic
representation dual to $\phi$'s). The notion of contragredient
does not really rear its head in the previous topics we have
covered, but is a feature of the more general functional
equations.

\subsection{Jacquet-Piatetski-Shapiro-Shalika (1979)}
\label{jpss}

Another proof of the analytic properties of the standard
$L$-functions of cusp forms on $GL(n)$ is a generalization of
Hecke's method (\secref{Hecke II}).   In the 1970s
Piatetski-Shapiro and Shalika \cite{PS,Shal} independently
developed their ``Whittaker'' expansions on $GL(n)$ in order to
generalize the expansion (\ref{FfromW}) of Jacquet-Langlands.  The
Whittaker function on $GL(n,\A_\Q)$ is given by the integral
\begin{equation}\label{glnwhit}
    W_\A\( g_\A \) \ \ = \ \ \int_{N(\Q)\backslash N(\A)} F_\A\(\(
    \begin{smallmatrix} 1 & u_{12} & u_{13} & \cdots & u_{1n} \\
                        {} &   1    & u_{23} & \cdots & u_{2n} \\
                        {}& & \ddots & \ddots &  \vdots \\
                        {} & &  & 1 & u_{\,n-1\,n} \\
                        {} & {} & {} & {} & 1
                   \end{smallmatrix}\) g_\A \)\,\overline{\psi(u_{12} + u_{23}
                   + u_{\, n-1\,n})}\,du \, ,
\end{equation}
with the integration over the subgroup $N$ of unit upper
triangular matrices.  The expansion of $F_\A$ in terms of $W_\A$
is given by
\begin{equation}\label{glnFfromW}
    F_\A(g_\A) \ \ = \ \ \sum_{\g\in N(\Q)\backslash P(\Q)}
    W_\A(\g\, g),
\end{equation}
where $P$ is the subgroup of $G$ consisting of matrices whose
bottom row is $(0 \ 0 \ \cdots \ 0 \ 1)$.  This also has an
explicit, classical version (``neoclassical,'' in the terminology
of Jacquet), which can be found in \cite{Bumpblue,jacquetindia} --
both of which are excellent references for this section.

 In a series of papers, Jacquet, Piatetski-Shapiro, and Shalika used
 these expansions to
generalize Hecke's construction to $L$-functions of automorphic
forms on $GL(n)$ (including properties {\bf E, BV}, and {\bf FE}),
and prove a converse theorem for $GL(3)$ (see
\cite{Bump-GL3,Bumpblue,jacquetindia,JPSS,JPSSrs}, and
\cite{ms-voronoi,ms-expos} for a different treatment).  This is a
big advantage over the method of Godement-Jacquet, whose integral
is over the large group $GL(n)$.  The integrals of Jacquet,
Piatetski-Shapiro, and Shalika instead involve integration over
one dimensional subgroups, matching the one complex variable of
the $L$-functions $L(s)$.   In later papers of Cogdell and
Piatetski-Shapiro, a powerful converse theorem has been
established  for $GL(n)$ (see \cite{C-PS} and \secref{la}). These
techniques lie close to the heart of the ``Rankin-Selberg''
method, which uses integral representations to generate a wide
variety of the Langlands $L$-functions we will come to in Sections
\ref{LaSh} and \ref{la} (see \cite{Bumpblue} for a thorough,
though slightly out of date, survey). While the statement of the
converse theorem is quite technical, it is similar in form to
Weil's \thmref{weilthm}, in that it involves the assumptions of
{\bf E}ntirety, {\bf B}oundedness in {\bf V}ertical strips, and
{\bf F}unctional {\bf E}quation; as in the above proof of Hecke's
theorem \thmref{heckecon}, these are used to shift a contour
integral which reconstructs an automorphic form using Mellin
inversion. However, an important difference is that their converse
theorem typically involves twisting by automorphic forms on
$GL(m)$, not merely Dirichlet characters (which, we have seen,
correspond to automorphic forms on $GL(1)$).  This is an important
topic in \secref{la}, where we state a typical version in
\thmref{ncrossn-2}; a general account can be found in
\cite{Cogdell}.

\section{Langlands-Shahidi (1967-)}\label{LaSh}

This section is meant for readers having some familiarity with Lie
groups, but it can be skipped without loss of continuity.
References include  \cite{ShahKore,Gel-Sha
book,HC,Borel,Moe-Wal,BernsteinPCMI,ShahidiPCMI}.  Our purpose
here is to describe the general method of obtaining analytic
properties of $L$-functions from Eisenstein series, generalizing
Selberg's method for $\zeta(s)$ in \secref{Sel}.  Many of the
applications in \secref{la} are based upon properties yielded by
the Langlands-Shahidi method.

  The theory of
Eisenstein series was widened by Langlands to more general Lie
groups in \cite{La1,La2}; in particular Langlands proved the
analytic continuation and functional equations that were useful in
Selberg's proof of the analytic properties of $\zeta(s)$. In his
Yale monograph \cite{Yale}, Langlands considered the constant
terms  of the completely general Eisenstein series. This time, a
wide variety of (generalized) $L$-functions appeared; his analysis
gives their meromorphic continuation. The calculations involved
are quite complicated and are performed adelically; they led
Langlands to define the $L$-group and ultimately to the
formulation of his functoriality conjectures.

Recall the example  of $GL(2)$ from \secref{Sel} (which we will
reconsider through group representations and adeles in
\secref{lashgl2ex}). There, analysis of the constant term and
first Fourier coefficient already sufficed for the analytic
continuation and functional equation of $\zeta(s)$ via Selberg's
method. Langlands proposed studying the non-trivial Fourier
coefficients in general, and Shahidi has now worked that theory
out
(\cite{shaduke,MR89h:11021,shannals,shahidi-ajm,ShahST,ShahKore,ShahidiPCMI,shahidi-icm})
along with Kim and others. In general it has been a difficult
challenge to prove the $L$-functions arising in the constant terms
and Fourier coefficients are entire. The analytic continuation of
Eisenstein series typically gives the meromorphic continuation to
$\C$, except for a finite number of poles on the real axis between
$0$ and $1$; these come from points where the Eisenstein series
themselves are not known to be holomorphic. A recent breakthrough
came with a clever observation of  H. Kim: the residues of the
Eisenstein series at these potential singularities are $L^2$,
non-cuspidal automorphic forms, and -- as in \secref{glnautsec} --
give rise to unitary representations that can be explicitly
described by the Eisenstein series they came from. Kim remarked
that results about the classification of irreducible unitary
representations show that many of these potential representations
do not exist, thus allowing one to conclude the holomorphy of the
Eisenstein series at these points in question! When combined with
\cite{shaduke,shannals}, this has recently led to new examples of
entire $L$-functions (more on this in \secref{la}).

\subsection{An Outline of the Method}

The following is a brief sketch of the main points of the method;
a fuller introduction with more definitions and detailed examples
can be found in \cite{ShahKore}.  Detailed examples of constant
term calculations can be found in many places, e.g.
\cite{Gel-Shabook,Yale,Blue,mil}. Though it is possible to
describe the method without adeles (as was done in \secref{Sel}),
their use is key in higher rank for factoring infinite sums and
product expansions into $L$-functions.  Because the
Langlands-Shahidi method utilizes various algebraic groups, we
will have to assume some familiarity with the basic concepts.  For
this reason we include an example of the $GL(2)$ case in
\secref{lashgl2ex}.

Let $F$ be a global field, $\A=\A_F$ its ring of adeles, and $G$ a
split algebraic group over $F$.  Much carries over to quasi-split
case as well, and we will highlight the technical changes needed
for this at the end.  Fix a Borel (= a maximal connected solvable)
subgroup $B\subset G$, and a standard maximal parabolic $P\supset
B$ defined over $F$.\footnote{The theory has an extension to
non-maximal parabolic subgroups, but this does not yield any extra
information about $L$-functions.  This matches the fact that the
Eisenstein series for maximal parabolic subgroups depend on one
complex variable, as do $L$-functions.} Decompose $B=TU,$ where
$T$ is a maximal torus. The parabolic can be also decomposed as
$P=MN$, where the unipotent radical $N\subset U$, and
 $M$ is the unique Levi component containing $T$.  Denote by
$^LG,^LM,^LN,$ etc. the Langlands dual $L$-groups (see
\cite{ShahKore} for definitions).

One of the key aspects of this method is that it uses many
possibilities of parabolics of different groups $G$, especially
exceptional groups.  This is simultaneously a strength (in that
there is a wide range of exotic possibilities) and a limitation
(in that there are only finitely many exceptional groups).

\subsubsection{Cuspidal Eisenstein Series}

Recall that an automorphic form in $L^2(\Gamma\backslash G)$ is
associated to a (unitary) automorphic representation of $G$.  Let
$\pi=\otimes_v\pi_v$ be a cuspidal automorphic representation of
$M(\A)$; we may assume that almost all components $\pi_v$ are
spherical unitary representations (meaning that they have a vector
fixed by $G(O_v)$, where $O_v$ is the ring of integers of the
local field $F_v$).  For these places $v$ the equivalence class of
the unitary representation $\pi_v$ is determined by a semisimple
conjugacy class $t_v \in \, ^LG$, the $L$-group.  This conjugacy
class is used to define the $L$-functions below in
(\ref{lfuncdef}).  The finite number of exceptional places are
where $\pi$ {\em ramifies}.

A maximal parabolic subgroup $P$ has a modulus character $\d_P$,
which is the ratio of the Haar measures on $M\cdot N$ and $N\cdot
M$. It is related to the simple root of $G$ which does not
identically vanish on $P$.  For any automorphic form $\phi$ in the
representation space of $\pi$, we can define the Eisenstein series

\begin{equation}\label{eisdef}
  E(s,g,\phi) \ \ = \ \ \sum_{\g\in P(F)\backslash G(F)}
  \phi(\g g)\,\d_P(\g
  g)^s
\end{equation}
and their constant terms
\begin{equation}\label{constdef}
 c(s,g,
 \phi) \ \ = \ \   \int_{N'(F)\backslash N'(\A)} E(s,ng,\phi)\, dn
 \, ,
\end{equation}
where $N'$ is the unipotent radical of the {\em opposite}
parabolic $P'$ to $P$ (it is related by the longest element in the
Weyl group). One can view the constant term as an automorphic form
on $M$, and we will shortly relate it to $\phi$ and $\pi$.  The
measure $dn$ is normalized to give the quotient $N'(F)\backslash
N'(\A)$ volume 1. The notion of constant term applies to any
parabolic, but $P'$ is the most useful one for our purposes.

\subsubsection{Langlands $L$-functions}
\label{rhoborn}

If $\rho$ is a finite-dimensional complex representation of $^LM$,
and $S$ is a finite set including the archimedean and ramified
places of $F$ and $\pi$, then the partial Langlands $L$-function
is

\begin{equation}\label{lfuncdef}
  L_S(s,\pi,\rho) \ \ = \ \ \prod_{v \notin S} \,
  \det(I\, - \, \rho(t_v)\,q_v^{-s})\i \, .
\end{equation}
Here $q_v$ is the cardinality of the residue field of $F_v$, a
prime power.  The full, {\em completed}, $L$-function involves
extra factors for the places in $S$, whose definition is technical
and in general difficult. This is connected to the local Langlands
correspondence, proven recently by Harris and Taylor for $GL(n)$
and by Jiang and Soudry for $SO(2n+1)$ (see
\cite{harris-taylor,MR1947454,jiang-soudry,langlandsreal,harris-icm,MR2002f:11178,MR2001i:11136,MR2001e:11052}).
 When $\rho$ is the standard
representation of $^LGL(n)=GL(n)$ and $F=\Q$, the Euler factors in
(\ref{lfuncdef}) agree with those in (\ref{degreen}); in general
the degree of $L_S(s,\pi,\rho)$ equals the dimension of $\rho$.

\subsubsection{The Constant Term Formula}

The constant term formula involves the sum of two terms.  The
first, which only occurs when the parabolic $P$ is its own
opposite $P'$, is $\phi(g)\d_P(g)^s$ -- simply the term in
(\ref{eisdef}) for  $\g = $ the identity matrix. Langlands showed
that the map from $\phi$ to the second term is described by an
operator
\begin{equation}\label{consttermformula}
M(s,\pi) \ \ = \ \ \(\prod_{j=1}^m \f{L(a_j\, s\, , \, \tilde{\pi}
\, , \, r_j)}{L(1\, + \, a_j \, s \, , \, \tilde{\pi} \, , \,
r_j)}\)\otimes_{v\in S} A(s,\pi_v) \, ,\end{equation} where the
$A(s,\pi_v)$ are a finite collection of operators, $r$ the adjoint
action of $^LM$ on the lie algebra of $^LN$, $r_1,\ldots,r_m$ the
irreducible representations it decomposes into, and $a_j$ integers
which are multiples of each other (coming from roots related to
the $r_j$).  The variety of decompositions of $r$ is what gives
this method much of its power for treating complicated
$L$-functions.  See \cite{ShahKore} for a fuller discussion, along
with an example for the Lie group $G_2$ and the symmetric cube
$L$-function. Tables listing Lie groups and the representations
$r_j$ occurring for them can be found in \cite{Yale} and
\cite{MR89h:11021}, for example.

\subsubsection{The Non-Constant Term: Local Coefficients}

We must now make a further restriction on the choice of $\pi$
involved, namely that it be {\em generic}, i.e. have a Whittaker
model.  This means that if $\psi$ is a generic unitary character
of $U(F)\backslash U(\A)$, we need to require

$$W(g,\psi) \ \ = \ \ \int_{U_M(F)\backslash U_M(\A)}
\phi(ng)\ \overline{\psi(n)} \ dn  \ \  \neq \ \ 0 \   ,  \ \ \ \
U_M \ = \ U \cap M$$ for some $\phi$ and $g$ (we have already seen
this notion in (\ref{adelwhitdef}) and (\ref{glnwhit})).

Shahidi's formula uses the Casselman-Shalika formula for Whittaker
functions (see \cite{cass-shal,shintani})  to express the
following non-constant term at the identity $g=e$ as

\begin{equation}\label{nonconstform}
  \int_{N'(F)\backslash N'(\A)}E(s,ne,\phi)\, \overline{\psi(n)}\,dn  \ \ =
  \ \
  \prod_{j=1}^m\f{1}{L(1\,+\,a_j\,s\, ,\, \tilde{\pi}\, , \, r_j)}
  \,
  \cdot  \, \prod_{v\in
  S}W_v(e) \, ,
\end{equation}
for a certain choice of $\phi$.  Applying the functional equation
of the Eisenstein series (which has the constant-term ratio
involved), one gets the ``crude'' functional equation for the
product of $m$ $L$-functions

\begin{equation}\label{bigprod}
  \prod_{j=1}^m L_S(a_j \,
  s \, , \, \tilde{\pi} \, ,\, r_j) \ \ = \ \
  \prod_{j=1}^mL_S(1\,-\,a_j\, s\, ,\, \pi\, ,\, r_j)
  \,  \cdot  \, \prod_{v\in S}\(\mbox{local
  factors}\)\,.
\end{equation}

Shahidi's papers \cite{shaduke} and \cite{shannals} match all the
local factors above to the desired $L$-functions (cf. the remark
after (\ref{lfuncdef})). This gives the full functional equation
for these $m$ $L$-functions, but only when multiplied together.
His 1990 paper \cite{shannals} uses an induction argument to
isolate the functional equation of each of the above $m$ factors
separately.

\subsubsection{Analytic Properties and the Quasi-Split Case}

It still remains to prove that the $L$-functions are entire,
except perhaps at $s=0$ and 1 (where the order of the poles is
understood, like for $\xi(s)$). The theory of Eisenstein series
provides this full analyticity for the $L$-functions arising in
the constant term  unless $\pi$ satisfies a self-duality
condition; even in this case, it can be shown that the
$L$-functions have only a finite number of poles, all lying on the
real axis between $0$ and 1.
 Kim's observation of using the unitary dual has
worked in many cases to eliminate this possibility. It is also
always possible to remove the potential poles by twisting by a
highly-ramified $GL(1)$ character of $\A_F$; this has been crucial
for applications to functoriality through the converse theorem
\cite{C-PS,MR95m:22009,MR97d:22019}, which we come to in
\secref{la}.

The main difference in the quasi-split case is that the action of
the Galois group $G_F$ is no longer trivial.   The $L$-groups are
potentially disconnected, inasmuch as they are  semi-direct
products of a connected component with $G_F$.  Also, the
representation $\rho$ used to define the Langlands $L$-functions
in (\ref{lfuncdef}) may also depend on the place $v$.

\subsection{$GL(2)$ Example}\label{lashgl2ex}

Here we reconsider the $\zeta$-function example from \secref{Sel},
but  in the framework of the Langlands-Shahidi method. In this
setting, the Eisenstein series on $G=GL(2)$ is defined by
\begin{equation}\label{eisgldef}
E(s,g,f) \ \ = \ \ \sum_{\gamma \in {B({\Q}) \backslash G({\Q})}}\
f({\gamma}g) \, ,
\end{equation}
where
 $M=P=P'=B= \left\{ \left(\begin{array}{cc} a& x\\
{}& b\end{array} \right) \right\} \subset G$ is the Borel
subgroup/minimal parabolic ($GL(2)$ is too small a group to afford
other interesting choices). The Eisenstein series formed from
$\pi$ are related to the representations {\em induced} from $\pi$,
from $M(\A)$ to $G(\A)$. In (\ref{eisgldef}) we may absorb the
factor $\d_P^{\ s}$ into $f$ by taking a vector in the {\em
induced} representation $I(s)=Ind_{B({\A})}^{G(\A)} \ |a|^s $,
which roughly speaking is the space of functions
$$
\left\{ f:G(\A) \rightarrow \C   \, \left|  \, f\(\ttwo ax{}{a\i}
g\)=|a|^{s+1}f(g) \right. \right\}.
$$
The Eisenstein series (\ref{eisgldef}) converges for $\Re\!\!(s)$
sufficiently large.  In fact, we may choose our $f\in I(s)$ so
that $E(s,g,f)$ reduces to just the classical Eisenstein series
$E(z,\f{1+s}{2})$ considered in \secref{Sel}.  To do this, we take
$f$ to be identically 1 on $\widehat{K}=O(2,\R)\times
\prod_{p<\infty} GL(2,\Z_p)$, and use the fact $B(\Q)\backslash
G(\Q) \simeq B(\Z) \backslash G(\Z)$\footnote{Classically
speaking, this isomorphism comes from the decomposition of any
rational matrix $g\in GL(2,\Q)$ as $g=bu$, $b\in B(\Q)$, $u\in
GL(2,\Z)$. The resulting indexing of $B(\Z)\backslash G(\Z)$ via
rational matrices gives powerful insight into how to arrange the
summands of the Eisenstein series into a computationally-useful
form.}; this, in view of the Iwasawa decomposition
$G=B\widehat{K}$, is the simplest choice. It corresponds to
(\ref{eisgldef}), taking $\pi$ to be the trivial representation of
$M(\A)$. In general, Eisenstein series are always induced from
automorphic forms on smaller groups, which in the example here is
just the constant function on the factors $M(\A)=GL(1,\A)\times
GL(1,\A)\subset B$.

To compute the constant term, we appeal to the Bruhat
decomposition
\begin{equation}\label{bru}
\gathered
  G \ \ = \ \ B \ \sqcup \ BwB \ \  = \ \  B
  \ \sqcup\  B wN\, , \\ w=\ttwo{}{-1}{1}{}\ \  , \ \ \  \ \
  N=\left\{\ttwo 1\star{}1     \right\}\subset G
\endgathered
\end{equation}
which is valid over any field.  When applied to $\g \in B(\Q)
\backslash G(\Q)$, it allows us to compute the constant term
integral over $N$ :
\begin{equation}\label{gl2intoverNbru}
  \aligned
  c(s,g,f) \ \ = \ \
\int_{N(\Q)\backslash N(\A)} E(s,ng,f)\, dn \ \ &=&f(g)& \ + \
\sum_{\g\in N(\Q)} \int_{N(\Q)\backslash N(\A)} f(w\, \g \,  n \,
g) \,dn
\\
&=&f(g)& \ + \ [M(s)f](g) \,,
  \endaligned
\end{equation}
 where
$M(s)$ is the {\em intertwining operator}
\begin{equation}\label{intertwine}
M(s)f(g) \ \ = \ \  \int_{N(\A)} \ f(w \, n \, g) \ dn
\end{equation} from $I(s)$ to $I(-s)$.
If $f$ is chosen to be a product $f(g)=\prod_p f_p(g_p)$, then the
integral (\ref{intertwine}) factors further to give an Euler
product, in analogy to (\ref{tatefold}) and (\ref{tateprod}).  For
the details of this example, see Langlands' article \cite{Blue};
in general, his constant term method \cite{Yale} gives similar
integrals for constant terms over general groups.  In any event,
for our example here where $f$ is trivial on $\widehat{K}$ (i.e.
so that $E(s,f,g)$ recovers the classical Eisenstein series),
$[M(s)f](e)=\f{\xi(s)}{\xi(s+1)}$.

 To complete our discussion let us  again compute  the $\psi$-th
Fourier coefficient of $E(s,g,f)$, where $\psi$ is a non-trivial
additive character of $N(\A)$ trivial on $N(\Q)$ (or,
equivalently, a non-trivial additive character of $\A$ trivial on
$\Q$). Afterwards we will use the functional equation of
$E(s,g,f)$ to get the functional equation of $\zeta(s)$.

We find
\begin{equation}\label{psithcoeffofE}
    \aligned
E_{\psi}(s,e,f) \ \ &:= \ \ \int_{N({\Q}) \backslash N(\A)}\
E(s,n,f) \ \overline{{\psi(n)}} \ dn  \\
&= \ \ {\frac{c(s)}{\xi(s+1)}}\, ,
    \endaligned
\end{equation}
and
$$
    E_{\psi}( \, -s \, , \, e \, , \, M(s)f \, ) \ \ = \ \
     \frac{\xi(s)}{\xi(s+1)} \ \frac{c(-s)}{\xi(1-s)}\, ,
$$
where $c(s)=c(-s)$ is non-zero (it is related to the $K$-Bessel
function appearing in (\ref{a1})). So using the functional
equation

$$E(s,\, e,\, f) \ \ = \ \ E(-s,\, e\, , M(s)f)\, ,$$
it follows that
  $$\xi(\,s\,) \ \ = \ \ \xi(\,1-s\,)\, .$$

\subsection{Boundedness in Vertical Strips and Non-vanishing}
\label{eisbvsec}

We have just seen how the functional equation and several analytic
properties of $L$-functions can be obtained via Langlands'
analytic continuation of Eisenstein series, which itself relies
 on spectral theory.  Through (\ref{psithcoeffofE}) and the
known holomorphy of $E(s,f,g)$ on the line $\Re{s}=0$ (which
follows from the general spectral analysis), one also obtains a
new proof of the famous result that $\zeta(s)$ never vanishes
along the line $\Re{s}=1$ (see
\cite{Jacquet-vanishing,sar-shalfest,moreno,shahidi-ajm}). Among
other things, this result is the key to the standard proof of the
Prime Number Theorem -- as was originally outlined by Riemann
himself in \cite{Riem}! In fact, this proof  of the non-vanishing
of $L$-functions on the line $\Re{s}=1$ using Eisenstein series
turns out to be the most general method available at present, in
some cases working far inside the known range of absolute
convergence of certain $L$-functions.

Intriguingly, it is possible to prove {\bf B}oundedness in {\bf
V}ertical strips using the Lang\-lands-Shahidi method.  This is
striking, because our other examples (Riemann, Hecke, and Tate,
e.g. \thmref{zetaebvfe}) all acquire {\bf BV} through an integral
representation of an $L$-function; here the $L$-function's
analytic properties are obtained very indirectly. Not
surprisingly, the argument is more round-about and subtle, but
pays off in that it turns out -- again -- to extend to more
general $L$-functions than treated by other methods alone (see
\cite{Gel-Sha bdness,sar-shalfest}).  This has been very important
for applications to the Langlands functoriality conjectures
through the converse theorem (see \secref{la} for further
discussion).

Recall that our  {\bf BV} assertion is that $s(s-1)\xi(s)$ is
bounded for $s$ in any vertical strip.
 In our situation, the
\emph{fact} that
      $$r(s)\, = \, \frac{\xi(s)}{\xi(s+1)}$$
satisfies the finite order estimate $O(e^{|s|^\rho})$ in $\Re{s}
\geq \frac{1}{2}$ is possible to prove using spectral theory. We
recall that $\xi(s)$ satisfies that finite order inequality
 in the
region $\Re{s}\ge 3/2$; this is because $|\zeta(s)|\le
\sum_{n=1}^\infty n^{-3/2} =\zeta(3/2)$ is bounded there, and
Stirling's formula (\ref{bigstir}) shows that
$\G(s/2)=O(e^{|s|^{1+\e}})$ there. Thus $\xi(s)=r(s)\xi(s+1)$
obeys the finite order inequality in $\Re{s}\ge 1/2$, and the
functional equation shows this is true for $\Re{s}\le 1/2$ as well
-- giving a new proof that $\xi(s)$ is of finite order. For
$\zeta$ and the general $L$-functions considered in \cite{Gel-Sha
bdness}, additional Eisenstein series and group representations,
along with some results of \cite{HC} and \cite{Mu}, are required.
In particular, more is needed than mere meromorphicity of the
Eisenstein series alone. One also needs an estimate on the growth
rate of Eisenstein series in vertical strips  going beyond the
original work of Langlands, which here is not being used  as
simply a ``black box."
 For a different, though related,
method of obtaining {\bf BV} and the non-vanishing of $\xi(s)$ on
the line $\Re\!\!(s)=1$ through Eisenstein series, see
\cite{Sarlet}.

\section{The Langlands Program (1970-)} \label{la}

Many articles have been addressed to the ``Langlands program''
(e.g. \cite{arthur,Gel-langlands survey} and the references in the
introduction), and it is not our desire to add to these. However,
one part of the program is closely connected to our discussion:
namely, it ``explains'' why the analytic continuation and
functional equation of the $L$-functions of automorphic forms on
$GL(n)$ \emph{probably} suffice  to ensure that \emph{any}
$L$-series in arithmetic has an analytic continuation and
functional equation!

\subsection{The Converse Theorem of Cogdell-Piatetski-Shapiro
(1999)}\label{cps1999}

 Let's see why: many arithmetic objects, such as elliptic curves,
 have an $L$-series attached to them  which are conjectured to be
 entire,
 and have functional equations similar to those possessed by
  our $L$-functions.  Some very interesting examples, which we
  will not touch on here directly, involve the Artin conjecture
  (which involves $L$-series of Galois representations); see
\cite{MR2001j:11026,MR98j:11106}.

   For an elliptic curve $E$
 defined over the rational numbers, this
 $L$-series, called the {\em Hasse-Weil} $L$-function $L(s,E)$,
  is defined by counting
 points on $E$ over varying finite fields.  Here $1+p-a_p$ is the number
 of points on the reduced curve modulo $p$, and $L(s,E)$ is
 defined by the Euler
 product in (\ref{eulerholom}) with $k=2$, except for a finite
 number of exceptional prime factors
 (see \cite{Silverman}).  The resemblance to the $L$-functions of
 holomorphic modular forms of weight 2 is the springboard for
   the celebrated ``Modularity Conjecture'' of Taniyama, Shimura,
  and Weil, which was  recently proven in
  \cite{Wiles}, \cite{Taylor-Wiles}, and \cite{BCDT}.  It
   asserts that
$L(s,E)=L(s,f)$, the $L$-function of some holomorphic cusp form
$f$ of weight 2 on
  $\G_0(N)$, where $N$ is a subtle invariant (the ``conductor'')
  calculable from the arithmetic of the curve.\footnote{Weil's
  contribution \cite{Weil} to the conjecture
   is closely related to \thmref{weilthm}: in fact his prediction
   of modularity on $\G_0(N)$, $N$ being the conductor of $E$,
   comes from a comparison of the expected functional equations of
   Hasse-Weil $L$-functions with the exact form of
   (\ref{weilcondiii}) in \thmref{weilthm}.
   }  Since the $L$-functions
 of modular forms
   are known to be entire through Hecke's theory (\secref{Hecke II}),
  we therefore now
  know that the Hasse-Weil $L$-functions of rational elliptic
  curves are entire.  Among other things, this gives a definition
 of $L(s,E)$ at the center of its critical strip, where the
  Birch-Swinnerton-Dyer conjecture
  asserts deep relations with arithmetic (\cite{BSD,gross-zagier,kolylog}).

One might ask if the modularity of an elliptic curve might be
proved using Weil's Converse \thmref{weilthm}.  Unfortunately,
this route requires one to know that the Hasse-Weil $L$-functions
are entire beforehand, which at present seems far beyond reach.
Thus the prospect of proving entirety and applying the Converse
Theorem to these arithmetic $L$-functions seems to be begging the
question. However, it is an interesting fact the Converse Theorem
on $GL(3)$ (proven in 1979 by Jacquet, Piatetski-Shapiro, and
Shalika) \cite{JPSS} {\em does} enter into the proof of the
Modularity Conjecture, as it had been earlier used to establish a
key step: the Langlands-Tunnell Theorem \cite{LT,Tun}.

Furthermore the Converse Theorem (in the form developed by
Piatetski-Shapiro and Cogdell, e.g. \cite{C-PS},
\thmref{ncrossn-2}) has had remarkable success towards the
Langlands Program in a different aspect. Roughly speaking, the
Langlands conjectures assert   correspondences between automorphic
forms on different groups. When starting with an automorphic form,
it is often possible to prove the {\bf E}ntirety, {\bf
B}oundedness in {\bf V}ertical strips, and {\bf F}unctional {\bf
E}quations of the (twisted) $L$-functions the Converse Theorem
requires. This has led to very significant progress, especially on
liftings of automorphic forms on $GL(n)$ to $GL(m)$, $n<m$.  In
particular, the recent breakthroughs of Kim and Shahidi (cf.
\cite{Kim-Sha3,KS-CR,Kim-Sha4,Kim-Sha,C-K-PS-S} and
\secref{LaSh}), and also of Lafforgue
(\cite{lafforgue,edfrenkel,laumon}), have proven many important
new examples of Langlands ``Functoriality'', by showing certain
$L$-functions are entire and then appealing to the Converse
Theorem of Cogdell and Piatetski-Shapiro (\cite{C-PS}).
Unfortunately the precise statements connected to the use of the
Converse Theorem are quite technical and complicated, and so we
will just make do with a less-technical (but still very useful)
version of the Converse Theorem in \thmref{ncrossn-2}. We will
then summarize the main applications in \secref{recent}.  In short
the basic idea, which can also been see through the standard
functoriality conjectures of Langlands, is the following: any
Langlands $L$-function (\ref{lfuncdef}) -- of any automorphic
form, on any group, over any field -- should itself be the
$L$-function of an automorphic form on $GL(n,\A_\Q)$.  So
$GL(n,\A_\Q)$ is speculated to the mother of all automorphic
forms, and its offspring $L$-functions are already known to have
an  analytic continuation and functional equation.

\subsection{Examples of Langlands $L$-functions: Symmetric Powers}
\label{examplesectionoflanglandslfunctions}

To get the full statements of the Langlands conjectures, one
requires even more than the definitions of
 the Langlands $L$-functions from (\ref{lfuncdef}).  For
 simplicity and the benefit of readers who have skipped
 \secref{LaSh}, we will explain the relevant $L$-functions here in
 the {\em everywhere-unramified} case, which corresponds to cusp
 forms on $GL(n,\R)$ invariant under $GL(n,\Z)$.

Recall the Euler product for the $L$-function of a modular form
$f$ from (\ref{firstlsfeul}); this formula is also valid for the
Maass forms from \secref{maass sec}, of course provided they are
eigenforms of all the Hecke operators $T_p$.  To better highlight
the symmetries involved (as well as to allow for more generality),
let us introduce the parameter $\b_p=\a_p\i$ and rewrite
(\ref{firstlsfeul}) as
\begin{equation}\label{secondlsfeul}
L(s,f)\ \ = \ \ \prod_{p} \, (1-\a_p \,p^{-s})\i \, (1-\b_p\,
p^{-s})\i.
\end{equation}
  Langlands has
defined higher degree Euler products  from  $L(s,f)$ called the
{\em symmetric $k$-th power} $L$-functions:

\begin{equation}\label{symk}
  L(s,Sym^k\,f) \ \ =\  \ \prod_p \prod_{j=0}^k (1- \, \a_p^{\,j}\
  \b_p^{\,k-j}\,
   p^{-s})\i.
\end{equation}

The definition of the symmetric power $L$-functions is a general
example of Langlands' method of creating new Dirichlet series out
of Euler products.  The major challenge, as we shall see, is to
derive important analytic properties of these  Dirichlet series,
and thereby put them on the same footing as the other
$L$-functions we have come across. His formulation is in terms of
finite dimensional representations, which in this case of $GL(2)$
is the $k+1$-dimensional representation on homogeneous polynomials
of degree $k$. Other examples give rise to Euler products which
are very symmetric, like this one on the righthand side of
(\ref{symk}).
 In
general, one starts by factoring the $L$-function of an
automorphic form on $GL(n)$ as
\begin{equation}\label{lspifact}
    L(s,\pi) \ \ = \ \ \prod_{p}\prod_{j=1}^n \,(1 \, - \,
    \a_{p,j}\,p^{-s})\i.
\end{equation}
New Euler products may be taken using symmetric combinations of
the $\a_{p,j}$ above (the individual $\a_{p,j}$ chiefly have
meaning only in the context of their aggregate $\{\a_{p,j}\}_{1\le
j \le n}$). In addition to the symmetric powers for $GL(2)$, there
are symmetric and exterior powers for $GL(n)$:
\begin{equation}\label{symgln}
    L(s,\pi,Sym^k) \ \ = \ \ \prod_p \prod_{i_1 \le i_2 \le \cdots \le
    i_k} (1 - \a_{p,i_1}\,\a_{p,i_2}\cdots \a_{p,i_k}\,p^{-s})\i
\end{equation}
\begin{equation}\label{extgln}
 L(s,\pi,Ext^k) \ \ = \ \ \prod_p \prod_{i_1 < i_2 < \cdots <
    i_k} (1 - \a_{p,i_1}\,\a_{p,i_2}\cdots
    \a_{p,i_k}\,p^{-s})\i\,.
\end{equation}
Given another $L$-function on $GL(m)$
\begin{equation}\label{lspiprime}
    L(s,\pi') \ \ = \ \ \prod_p \prod_{k=1}^m\,(1 \, - \,
    \b_{p,k}\,p^{-s})\i\, ,
\end{equation}
Langlands forms the ``Rankin-Selberg'' tensor product
\begin{equation}\label{lspipi'}
    L(s,\pi\otimes \pi') \ \ = \ \ \prod_p \prod_{j=1}^n
    \prod_{k=1}^m \, (1 \, - \, \a_{p,j}\,\b_{p,k}\,p^{-s})\i\, ,
\end{equation}
in analogy with the classical constructions
\cite{rankin,selberg-rs} for $GL(2)$ (see \cite{Bumpblue}). There
is a complementary theory for $\G$-factors and completed, global
Langlands $L$-functions as well.   The general Langlands
construction is in terms of finite dimensional representations of
$L$-groups (\secref{rhoborn}); in particular, they can be repeated
in various configurations.  Now, thanks to the recent proof of the
local Langlands correspondence by Harris and Taylor for $GL(n)$
\cite{harris-taylor,MR1947454,langlandsreal,harris-icm,MR2002f:11178,MR2001i:11136,MR2001e:11052},
the definitions at the ramified places can be made also.
Langlands' deep conjectures, in these cases here, assert that each
of the $L$-functions defined above is in fact the $L$-function of
some automorphic form on some $GL(d)$, where $d$ is the degree of
the Euler product in each case (i.e. the number of factors
occurring for each prime).  Or, in other words, if his
symmetric-looking Euler products look like the Euler product of an
automorphic form as in (\ref{lspifact}), they probably are!

\subsection{Recent Examples of Langlands Functoriality (2000-)}
\label{recent}

To continue, we now wish to focus on the examples of Langlands'
lifting mentioned above.  We will describe various lifts which
start with automorphic forms on $GL(n)$, and create automorphic
forms on some $GL(m)$, $m>n$.  Though many examples of Langlands
functoriality are known in various types of cases, this class is
very analytic and has largely been unapproachable without using
the types of analytic properties of $L$-functions that we have
come across in this paper. When considered as correspondences
between eigenfunctions one space and another, the lifts below are
quite stunning theorems in harmonic analysis,  made possible by a
deep use of the arithmetic of $GL(n,\Z)$.

Having explained the tensor product $L$-function (\ref{lspipi'}),
we can now state a version of the converse theorem (in practice,
slightly weaker assumptions are often used, as well as analogs
over different number fields):

\begin{thm}\label{ncrossn-2}($GL(n)\times GL(n-2)$ Converse
Theorem -- \cite{C-PS})

Consider the Euler product $L(s,\pi)$ (\ref{lspifact}), and assume
that it is convergent for $\Re\!{s}$ sufficiently large. Suppose
that $L(s,\pi)$ along with all possible tensor product
$L$-functions $L(s,\pi \otimes \tau)$, for $\tau$  an arbitrary
cuspidal automorphic form on $GL(m,\A_\Q)$,  $1\le m \le n-2$,
satisfy {\bf E}ntirety, {\bf B}oundedness in {\bf V}ertical
strips, and the {\bf F}unctional {\bf E}quation. Then $L(s,\pi)$
is in fact the $L$-function of a cuspidal automorphic form on
$GL(n,\A_\Q)$.
\end{thm}

Of course, in this statement we have not described the global
$L$-function (e.g. $\G$-factors) whose analytic properties we are
describing, but it is similar to the ones from \secref{Hecke II}.
Needless to say, \thmref{ncrossn-2} is a generalization of
\thmref{weilthm}.  When $n=3$, it is an earlier theorem of
Jacquet-Piatetski-Shapiro-Shalika \cite{JPSS}.  To use the
Converse Theorem to establish lifting to $GL(n)$, one still needs
to show that various tensor product $L$-functions obey the
analytic conditions it requires.  Such properties are themselves
very difficult assertions in their own right, and progress has
been hard won.  We shall now describe the established lifts from
$GL(n)$ to $GL(m)$ that were mentioned at the end of the last
subsection.

 The first such example is the symmetric square
lift $Sym^2: GL(2)\rightarrow GL(3)$, the so-called
Gelbart-Jacquet lift \cite{gel-jaclift}.   Because it is the
simplest of these to explain, we will spend a moment to go over
how it is proved.  A central role is played by Shimura's integral
representation of the symmetric square $L$-function
\cite{shimurasymsq}; one obtains the necessary analytic properties
of $L(s,Sym^2\pi \otimes \chi)$, where $\chi$ is a Dirichlet
character (recall that these are automorphic forms on $GL(1)$).
Then the Converse Theorem of \cite{JPSS} (i.e. \thmref{ncrossn-2}
with $n=3$) implies the existence of a cuspidal automorphic
representation $\Pi$ whose $L$-function $L(s,\Pi)=L(s,Sym^2\pi)$
-- i.e., the Langlands functorial symmetric square lift from
$GL(2)$ to $GL(3)$.

Examples on $GL(n)$ for $n\ge 4$ require the
Cogdell-Piatetski-Shapiro versions of the Converse theorem, and
are quite technical, even in description.  Fortunately, many have
been achieved in the last few years, mainly as a consequence of
new analytic properties from the Langlands-Shahidi method
(\secref{LaSh}), mined from various configurations of parabolic
subgroups in exceptional groups such as $E_8$.  Here is a summary
of the recent lifts:
\begin{thm}
\label{recent lifts}

The following instances of Langlands functoriality are known. That
is, in each case there are automorphic forms on the target $GL(n)$
whose standard $L$-functions agree with the  Langlands
$L$-functions on the source group  (cf.
\secref{examplesectionoflanglandslfunctions}):

\begin{itemize}

\item Gelbart-Jacquet \cite{gel-jaclift}.  $Sym^2: GL(2)
\rightarrow GL(3)$.

\item Ramakrishnan \cite{ramlift}.  Tensor Product: $GL(2) \times
GL(2) \rightarrow GL(4)$.

\item Kim-Shahidi \cite{Kim-Sha3,KS-CR}.   Tensor product:
$GL(2)\times GL(3) \rightarrow GL(6)$.

\item  Kim-Shahidi \cite{Kim-Sha3,KS-CR}.  $Sym^3: GL(2)
\rightarrow GL(4)$.

\item Kim \cite{KS-CR,Kim-Sha4}. $Ext^2: GL(4)\rightarrow GL(6)$
weakly automorphic (bad at 2 and 3).

\item Kim \cite{KS-CR,Kim-Sha4}. $Sym^4: GL(2)\rightarrow GL(5)$.

\item Cogdell-Kim-Piatetski-Shapiro-Shahidi
\cite{C-K-PS-S,C-K-PS-S2}: Weak functoriality to $GL(n)$ for
generic cusp forms on split classical groups.
\end{itemize}
\end{thm}

The notion of ``weak'' automorphy  means that an automorphic form
on the target $GL(n)$ exists whose  $L$-function matches the
desired Euler product  -- but {\em except} perhaps for a finite
number of factors.
 Much more about
these results can be found in these references, and also the ICM
lectures \cite{shahidi-icm,cogdell-icm}. Ramakrishnan's result
used an integral representation for a triple-product $L$-function
(\cite{garrett-annals-1987,ikeda,ps-rallis-triple,harris-kudla}),
but can also now be proven using the Langlands-Shahidi method
(\cite{Kim-Sha4}). The last example mentioned here is a lift from
 generic cuspidal automorphic forms on $SO(n)$ or $Sp(2n)$
to some  $GL(m)$ (see
\cite{cogdell-icm,shahidi-icm,C-K-PS-S,C-K-PS-S2}).  A differing
``descent method'' (i.e. studying the opposite direction of the
lift) of Ginzburg, Rallis, and Soudry \cite{MR2002g:11065} can be
used to establish  the lifts of \cite{C-K-PS-S,C-K-PS-S2} in
strong form; in other words, the adjective ``weak'' can be removed
from the last assertion of \thmref{recent lifts}.

Of course  Langlands' conjectures go far beyond these examples
involving only $GL(n)$ over a number field.  Other routes, through
theta liftings (see \cite{MR2003c:11051,MR92j:11045}) and the
Arthur-Selberg trace formula (see
\cite{arthur,MR98d:22017,MR98j:11105}), have also provided many
instances of Langlands Functoriality.  In particular, the trace
formula is in some sense the most successful when successful, in
that it gives a very complete description and characterization of
the lifts it treats. Nevertheless, the full force of Langlands'
Conjectures seem absolutely beyond current technology (see
\cite{beyondendoscopy} for intriguing comments by Langlands on the
limitations of the trace formula). We shall not describe these nor
the exact formulations of the Converse Theorem here, though we
hope we have transmitted the flavor of the arguments and technical
analytic properties such as {\bf EBV} which have put these recent
results in grasp.

\subsection{Applications to Number Theory (2001-)}

The coefficients of modular forms on the complex upper half plane
$\U$ play a fundamental role in many problems in number theory.
For example, the coefficients of holomorphic modular forms
$\phi(z)=\sum_{n\ge 0} c_n\,e^{\,2\,\pi\, i \, n \, z}$ can be
related to various counting problems, such as the number of ways
to represent a number as a sum of squares, or the number of points
on an elliptic curve (\secref{cps1999}).  The coefficients $a_n$
of the non-holomorphic Maass forms in (\ref{intmaass}) are also
related to number theory as well, ranging for example from Galois
theory to the properties of Kloosterman sums $\sum_{x\bar{x}\equiv
1\pmod p} e^{2\,\pi \, i \,\f{a x + b \bar{x}}{p}}$
\cite{sel1965,Sarbook,greeniwaniec,Iw-Sar,goldfeldsarnak}.  The
sizes of the $a_n$ and   eigenvalue parameter $\nu$, along with
their distributions, are very important in many instances; in the
remainder of this section, we will describe the role of the
analytic properties of $L$-functions in gleaning some of this
information.

\subsubsection{Progress towards the Ramanujan and Selberg
conjectures}

Recall Ramanujan's $\D$ form, defined in
 (\ref{ramformdef}).  We mentioned that Ramanujan conjectured a
 bound on the normalized coefficients $a_n$ of his $\D$ form, a
 bound which has a natural generalization to the coefficients of
 modular forms of any weight, and to Maass forms as well.

 \begin{conj}\label{ramanujan}(Ramanujan Conjecture)
Let $\phi(x+iy)$ be either a holomorphic cusp form of weight $k$
with Fourier coefficients $c_n = a_n n^{(k-1)/2}$ as in
(\ref{cnholom}), or a Maass form with Fourier coefficients $a_n$
as in (\ref{intmaass}). Then $a_n = O(|n|^\varepsilon)$ for any
$\e>0$ (of course the implied constant in the $O$-notation here
may depend on $\varepsilon$). When $\phi$ is a Hecke eigenform and
$a_1$ is normalized to be 1, then equivalently $|a_p|\le 2$.
 \end{conj}

This conjecture was proven in the holomorphic case by Deligne
\cite{Deligne} in 1974, but remains open for Maass forms.

Years later after Ramanujan, Selberg made a separate conjecture
about the size of the parameter $\nu$ that enters into the Fourier
expansion of Maass forms.  It is related to the Laplace eigenvalue
by $\l=1/4-\nu^2$. Selberg conjectured
\begin{conj}(Selberg, 1965 \cite{sel1965})\label{selconj}
Let $\l>0$ be the Laplace eigenvalue of a Maass form for
$\G\backslash \U$, where $\G$ is a congruence subgroup of
$SL(2,\Z)$.  Then $\l \ge \f 14$ (equivalently, $\nu$ is purely
imaginary).
\end{conj}
 Selberg was originally
motivated by questions involving cancellation in sums of
Kloosterman sums, but his question is a deep one about the nature
of the Riemann surfaces $\G\backslash \U$.  Their volumes grow to
infinity as their index increases, and one would naively expect an
accumulation of small Laplace eigenvalues.  However, Selberg
predicts a barrier at $\l=\f 14$.  This has implications for the
geometry of $\G\backslash \U$: intuitively, small eigenvalues are
a measure of how close a surface is to being disconnected, since,
after all, the multiplicity of the eigenvalue $\l=0$ is the number
of disconnected components.  These ideas have played a crucial
role in the development of {\it expander graphs}: discrete
combinatorial networks which have relatively few edges connecting
their vertices, but which are extremely difficult to disconnect by
removing only a moderate number of edges.  See
\cite{LPS,Margulis,Lbook,Sarbook,sarnotices,Murty}.
 By the way, Maass forms with eigenvalue exactly
equal to $\f 14$ are known to exist, and in fact they will play a
role later at the end of this section.  So Selberg's conjecture,
if true, is sharp!

Not long after Selberg's conjecture, Satake observed a unifying
reformulation of both the Ramanujan and Selberg conjectures, in
terms of representation theory (more specifically, tempered
representations).  Suppose $\phi$ is a Hecke eigenform.  A key
idea was the parametrization of the
 Hecke eigenvalues (which are also Fourier coefficients)
  $a_p$, for $p$ prime, as
$a_p=\a_p+\a_p\i$, $\a_p \in \C$ (cf. (\ref{firstlsfeul})). This
has  meaning from the representation theory of the group
$GL(2,\Q_p)$, and is analogous to the convention of writing the
Laplace eigenvalue $\l=\f 14 - \nu^2$. The statement that
$|a_p|\le 2$ is equivalent to proscribing that the complex modulus
satisfy $|\a_p| = 1$. Writing $\a_p$ as $p^{\mu_p}$, the
connection between Ramanujan's and Selberg's conjectures becomes
even more clear: both $\nu$ and all $\mu_p=\log_p (\a_p)$ should
be purely imaginary.

The generalized Ramanujan-Selberg conjecture asserts this
phenomenon holds for $GL(n)$:

\begin{conj}\label{genramconj}
If $\pi$ is a cusp form on $GL(n)$ which is unramified at the
prime $p$, the quantities $\a_{p,j}$ appearing in the Euler
product (\ref{lspifact}) obey $|\a_{p,j}| = 1$.
\end{conj}

All but a finite number of primes are ramified for $\pi$, and none
of them are when $\G=GL(n,\Z)$.  A similar statement for the
archimedean case $p=\infty$ is conjectured to be true for the
parameters $\mu_{\infty,j}$, generalizing $\nu$, that appear in
the $\G$-factors that multiply the $L$-function in its global,
completed form (see (\ref{glngam})).   For the cognoscenti we will
note that the Ramanujan conjecture~\ref{genramconj} has a
statement in terms of representation theory which covers the
ramified places as well: the associated local representations
$\pi_p$ of $GL(n,\Q_p)$ should all be {\em tempered}.

Though the generalized conjecture for $GL(n)$ is of course no
easier than it was for $GL(2)$, this added perspective has been
crucial for two reasons. The first is that we know a ``trivial''
or ``local bound'' coming from representation theory
\cite{jacquet-shalika} that $|\Re{\mu_{p,j}}| < 1/2$ for all
places $p\le \infty$.  For $GL(2)$ this is quite trivial indeed:
it states, for example, that the Laplace eigenvalue is merely
positive, and that the corresponding bound on the Hecke eigenvalue
$a_p$ comes directly from the boundedness of a cusp form. However,
for $GL(n)$, $n>2$, this bound actually becomes quite deep, due to
a separation feature between the trivial and non-trivial unitary
irreducible representations of $GL(n)$.

The second reason is that the Langlands program connects
automorphic forms on different $GL(n)$'s, for example through
symmetric powers. Notably, the factors $(1-\a_p^{n-1} p^{-s})\i$
and $(1-\a_p^{1-n}p^{-s})\i$ occur in the Euler product for the
$n-1$-st symmetric power from $GL(2)$ to $GL(n)$ (formula
(\ref{symk})). If this symmetric power $L$-function was indeed the
$L$-function of a cusp form, we would conclude that
$p^{-1/2}<|\a_p^{n-1}| < p^{1/2}$ from the ``trivial bounds''
above.  This gives an improved bound towards the Ramanujan-Selberg
conjectures for any $n$ for which the symmetric power lifting can
be established -- a bound which approaches the conjecture
$|\a_p|=1$  itself as $n\rightarrow \infty$.  A similar
magnification can be performed with the archimedean parameters
$\mu_{\infty,j}$.

Thus the Langlands program (in particular, the symmetric power
functorial liftings) implies both the Ramanujan and Selberg
conjectures, and their generalizations to $GL(n)$!  It should be
noted that this strategy is  different from Deligne's and other
arguments coming from algebraic geometry -- which themselves have
been successful for certain {\em cohomological} forms, but do not
apply to  Maass forms, for example.  (Actually Deligne's argument
uses the ``magnification'' mechanism of the previous paragraph,
but in a different context.)  Anytime a new lift is proven or a
new bound on the $|\a_{p,j}|$ of  cusp forms on $GL(n)$ is
established, it results in a bound towards the Ramanujan and
Selberg conjectures. Using results of \cite{LRS1}, an analytic
technique of \cite{Duke-Iwaniec}, and the recent progress of
Kim-Shahidi described in the previous subsection, the following
bounds have been proven:

\begin{thm}(Kim-Sarnak  \cite[Appendix 2]{Kim-Sha4})
If $\pi$ is a cusp form on $GL(2,\A_\Q)$ then
\begin{equation}\label{ks-localp}
    p^{-7/64} \ \le \ |\a_p| \ \le \ p^{\,7/64},  \ \ \ \text{~if $\pi$ is
    unramified at $p$,}
\end{equation}
and
\begin{equation}\label{ks-localinfty}
    \l \  \ge  \  \ \f{975}{4096} \ \  = \ \ \f 14 \
     -\  \(\f{7}{64}\)^2 \ \ \approx \ \  .238037 \ , \ \ \ \text{~if $\pi$
    comes from a Maass form}.
\end{equation}
\end{thm}

This theorem is for $\Q$, but results are also possible over
general number fields.  A weaker estimate (with $\f{7}{64}$
replaced by $\f 19$) is established by Kim-Shahidi in
\cite{Kim-Sha}, using their recent progress and techniques from
\cite{MR89h:11021}.

\subsubsection{The distribution of the Hecke eigenvalues,
Sato-Tate}

Having seen that the Hecke eigenvalue parameters $\a_p$ for a
$GL(2)$ modular, Hecke eigenform should lie on the unit circle in
the complex plane, we now turn to their distribution over this
circle as $p$ varies.  The question has its origin in conjectures
and investigations made independently by Sato and Tate
\cite{tatewoodshole} for the $a_p$ of rational elliptic curves
(\secref{cps1999}). Namely, if we consider the phase of $\a_p$,
i.e. the angle $\theta_p$ such that $a_p=2\cos\theta_p$, the
$\theta_p \in [0,\pi]$ should be equidistributed with respect to
the measure $\f {2}{\pi} \sin^2\theta \, d\theta$. This means that
\begin{equation}\label{satotate}
\lim_{X\rightarrow \infty} \  \f{\#\{  \, \a < \theta_p < \b \,
\mid  \, p\le X \} }{\#\{ \,  p\le X
    \}}      \ \ = \ \ \int_{\a}^\b \, \left[ \f{2}{\pi} \sin^2\theta
    \right]d\theta \,;
\end{equation}
when viewed in terms of the $a_p$ themselves, the conjecture
states that a histogram of the $a_p$ is governed by the
distribution $\f{1}{2\pi}\sqrt{4-x^2}$, which looks like a
semi-circle (really, semi-ellipse) between $-2$ and $2$ of area 1.
The Sato-Tate semi-circle measure occurs in many contexts; here it
is related to the Weyl integration formula, which weighs the
relative sizes of conjugacy classes in $SL(2,\R)$.  This
conjecture is not meant to be valid for {\em all} elliptic curves
(nor, by extension, to all modular forms via Wiles et al), but
instead only for the ``typical'' (i.e. non-CM) elliptic curve.  In
the other cases, the distribution is much simpler and the desired
results are known (see \cite{serrelad}).  Regardless, the
Sato-Tate conjecture is expected to also hold for most modular and
Maass forms, as we shall see shortly.  See \cite[p.210]{marchen}
for a generalization to $GL(n)$.

As in nearly all distributional questions in number theory, an
equivalent formulation of the Sato-Tate conjecture can be made in
terms of the {\em moments}
\begin{equation}\label{stmoment1}
   S_m(X) \ \ := \ \   \sum_{p\le X} a_p^m\,.
\end{equation}
Conjecturally $S_m(X)/\pi(X)$ should tend to the constant
\begin{equation}\label{stmoment2}
\lim_{X\rightarrow\infty} \ \f{S_m(X)}{\pi(X)} \ \
 = \ \ \f 1{2\,\pi}
\int_{-2}^2 x^m\,\sqrt{4 -x^2}\,dx  \end{equation} ($\pi(X)$, as
in the introduction, refers to the number of primes $p\le X$).  In
fact, the truth of (\ref{stmoment2}) for all $m\ge 0$ implies the
Sato-Tate conjecture (\ref{satotate}).

 In view of the connection with the powers of $a_p = \a_p+\a_p\i$,
  and symmetric power $L$-functions in the previous subsection, it
is not surprising that symmetric power $L$-functions play a role
in the Sato-Tate conjecture as well. In fact, the non-vanishing
and holomorphy of the $m$-th symmetric power $L$-function
$L(s,Sym^m\,\pi)$ in the region $\Re{s}\ge 1$ implies the $m$-th
moment (\ref{stmoment2}) (see \cite{serrelad,ShahST}, and also
Ogg's paper \cite{ogg2}, which shows the holomorphy actually
implies non-vanishing). Kim and Shahidi have now established this
for $m\le 9$ (\cite{Kim-Sha}). Actually the nonvanishing of
$\zeta(s)$ in the region $\Re{s}\ge 1$ is essentially what Riemann
observed implies the prime number theorem $\pi(X)=\sum_{p\le X}1
\sim X/\log X$, so it is natural to see this analytic condition
appear in a counting problem.  This is a typical way exotic
$L$-functions enter into analytic number theory.

Though the formulation of (\ref{satotate}) here implicitly assumed
the Ramanujan conjectures, (\ref{stmoment2}) is more general. It
can be viewed as saying that the Ramanujan conjecture is true on
average -- and much more.   We note in passing that various
on-average results can be proven using the theory of
$L$-functions. The Rankin-Selberg method \cite{rankin,selberg-rs}
has its origin in this issue for $GL(2)$; the generalization of
the Rankin-Selberg method to $GL(n)$ (\cite{JPSSrs,shaduke})  also
gives a weaker on-average version of Ramanujan.  One can also give
a relatively large lower bound on the percentage of primes $p$
such that the Ramanujan conjecture holds for $p$
(\cite{rama-lowerbd,Kim-Sha,Ramakrishnan}).

Finally, we conclude by describing a result of Sarnak
\cite{sarint}. We had mentioned before that the Maass forms for
$SL(2,\Z)\backslash \U$ -- the non-holomorphic
$L^2$-eigenfunctions of the Laplacian -- are quite mysterious in
nature, and none has been explicitly described.  However Maass, in
his original paper \cite{maass}, constructed some examples for
$\G\backslash \U$, where $\G$ is a congruence subgroup of
$SL(2,\Z)$.  These, and generalizations coming from Galois theory
through the Artin conjecture, are very special types of Maass
forms, and come from algebraic constructions. In particular,
several give fascinating examples of Maass forms (\ref{intmaass})
whose coefficients $a_n$ are relatively small integers -- bounded
in absolute value by the number of divisors of $n$.  This is
remarkable because of the discreteness and limitation of the
possible coefficients.  In \cite{sarint} Sarnak considers
hypothetical Maass forms with integral coefficients that are {\em
not} examples of the known constructions from Galois theory. In
these cases, the results of Kim and Shahidi \cite{Kim-Sha} on the
non-vanishing and holomorphy of $L(s,Sym^m \pi)$ on the line
$\Re{s} = 1$ give the asymptotics of the $m$-th moment, i.e.
 (\ref{stmoment2}), for $m\le 9$.
 If the
coefficients are indeed integral, the Ramanujan conjecture asserts
that the $a_p$ should only assume one of the five values
$\{-2,-1,0,1,2\}$. This constraint makes it difficult to match the
predicted moments, and in fact with $m=6$ it is possible to show
the impossibility of all the $a_p$ being integral.  Indeed, even
without the Ramanujan conjecture, the assumption that all
$a_p\in\Z$ can be ruled out simply by taking linear combinations
of (\ref{stmoment2}), and concluding that
\begin{equation}\label{}
    \lim_{X\rightarrow \infty} \,\f{1}{\pi(X)} \, \sum_{p\le X}
    P(a_p) \ \ =
    \ \ \f{1}{2\pi} \, \int_{-2}^2 P(x) \, \sqrt{4-x^2}\, dx  \ \
     = \ \ 1\,.
\end{equation}
Here $P(x)=x^2(4-x^2)(x^2-1)$, a sixth degree polynomial which
vanishes at the integers $\{-2,-1,0,1,2\}$, and is negative at all
others; a contradiction arises because the righthand side is
positive.  As a result, one obtains the first algebraicity result
in the subject of Maass forms:

\begin{thm}(Sarnak \cite{sarint})\label{sarnakstheorem}
    Let $\phi$ be a Maass form for $\G\backslash \U$ as in
    (\ref{intmaass}) with integral coefficients.  Assume $\phi$ is a
    Hecke eigenform. Then $\phi$
    arises from a Galois representation, and in particular the
    Laplace eigenvalue of $\phi$ is $\f 14$ (i.e. $\nu=0$).
\end{thm}

A generalization has been established by Brumley \cite{Brumley}.
We leave the reader with a some open conjectures -- both widely
believed to be true and supported by numerical evidence -- on
which the ideas of functoriality and $L$-functions have brought an
interesting perspective.

\begin{conj}
 (See \cite{cass-gln})
  Let $\phi$ be a Maass form which has
   Laplace eigenvalue $\f 14$.
Does $\phi$ necessarily arise from a Galois representation?
\end{conj}

\begin{conj}(Cartier, \cite{cartier-numeric}) Is the
Laplace spectrum of Maass forms for
 $SL(2,\Z)\backslash \U$ simple?  In other words, can there be
 two linearly-independent Maass forms on $SL(2,\Z)\backslash \U$
 sharing the same eigenvalue?
\end{conj}

\vspace{3 cm}

 \begin{tabular}{lcl}
Stephen S. Gelbart &&Stephen D. Miller                   \\
Faculty of Mathematics and Computer Science
&&Department of Mathematics    \\
Nicki and J. Ira Harris Professorial Chair && Hill Center-Busch Campus \\
The Weizmann Institute
of Science &&    Rutgers University          \\
Rehovot 76100&&        110 Frelinghuysen Rd             \\
Israel &&   Piscataway, NJ 08854-8019                  \\
{\tt gelbar@wisdom.weizmann.ac.il} &&   {\tt
miller@math.rutgers.edu}
\end{tabular}

\end{document}